\documentclass[letterpaper]{article}
\usepackage{amssymb}
\usepackage{amsmath}
\usepackage{amsfonts}
\usepackage{amscd}
\usepackage{latexsym}
\usepackage{graphicx}
\usepackage{amsthm,enumerate}
\usepackage[all]{xy}

\newtheorem{theorem}{Theorem}[section]
\newtheorem{proposition}[theorem]{Proposition}
\newtheorem{lemma}[theorem]{Lemma}
\newtheorem{corollary}[theorem]{Corollary}
\newtheorem{definition}{Definition}[section]

\newcounter{example}
\renewcommand{\theexample}{\arabic{example}}

\begin{document}
\title{Open Orbits and\\ Augmentations of Dynkin Diagrams}
\author{Sin Tsun Edward Fan and Naichung Conan Leung}
\date{5 March, 2009}
\maketitle
\begin{abstract}
Given any representation $V$ of a complex linear reductive Lie
group $G_{0}$, we show that a larger semi-simple Lie group $G$
with
\[
\mathfrak{g}=\mathfrak{g}_{0}\oplus V\oplus V^{\ast }\oplus \cdots
,
\]%
exists when $V$ has a finite number of $G_{0}$-orbits together
with a few exceptions corresponding to a twisted version of it. In
particular, $V$ admits an open $G_{0}$-orbit. Furthermore, this
corresponds to an augmentation of the Dynkin diagram of
$\mathfrak{g}_{0}$.

The representation theory of $\mathfrak{g}$ should be useful in
describing the geometry of manifolds with stable forms as studied
by Hitchin.
\end{abstract}

\section{Introduction} \label{sect:1}
Each fundamental representation $\Lambda ^{k}\mathbb{R}^{n}$ of
$GL_n( \mathbb{R})$ corresponds to the node labelled by $k$ in the
following Dynkin diagram $\Gamma $ of $GL_n( \mathbb{R})$
\begin{displaymath}
\xymatrix{
\circ\ar@{-}[r]_{\!\!\!\!\!\!\!\!\!\!\!\!\!\!\!\!\!\!\!1} &\circ \ar@{-}[r]_{\!\!\!\!\!\!\!\!\!\!\!\!\!\!\!\!\!\!\!\!2} &\cdots \ar@{-}[r] &\circ \ar@{-}[r]_{\!\!\!\!\!\!\!\!\!\!\!\!\!\!\!\!\!n-2} &\circ \ar@{-}[l]^{\,\,\,\,\,\qquad \qquad n-1}\\
}.
\end{displaymath}
It is interesting to observe that $\Lambda ^{k}\mathbb{R}^{n}$ has
an open orbit precisely when we can form a new Dynkin diagram by
attaching a new node to $\Gamma $ at the place labelled by $k$.
Furthermore, the simple Lie algebra $\mathfrak{g}$ corresponding
to this new
Dynkin diagram can be built from $\mathfrak{gl}_n$ and $\Lambda ^{k}\mathbb{R}^{n}$ and it is of the form%
\[
\mathfrak{g}=\mathfrak{gl}_n \oplus \Lambda
^{k}\mathbb{R}^{n}\oplus \left( \Lambda ^{k}\mathbb{R}^{n}\right)
^{\ast }\oplus \cdots .
\]

In this paper, we show that this phenomenon holds true in general.
Given any complex linear reductive Lie group $G_{0}$ and any
irreducible representation $V$ of it. One could try to form a
larger semi-simple Lie group $G$, or equivalently a Lie
algebra $\mathfrak{g}$, of the form%
\[
\mathfrak{g}=\mathfrak{g}_{0}\oplus V\oplus V^{\ast }\oplus \cdots
,
\]%
as a Lie algebra with a $\mathbb{Z}$-gradation.

The main result of this paper shows that such a Lie algebra
$\mathfrak{g}$ exists precisely when the number of $G_{0}$-orbits
in $V$ is finite. Moreover, the Dynkin diagram of $\mathfrak{g}$
is an augmentation of the Dynkin diagram of $\mathfrak{g}_{0}$ in
the same way as in the $GL_n$ case, or a twisted version of it.
Furthermore, the length of the $\mathbb{Z}$-gradation can be
easily read off from the Kac diagram \cite{ref12} and it is at
most six. In particular, $V$ has an open orbit. Irreducible
representations which admit open orbits are completely classified
(\cite{ref13}, also see Tables \ref{table:2} and \ref{table:3} for
those which admits a finite number of orbits). We see that all
cases except one have a finite number of orbits. More precisely,
there are one series of such representations, namely $(GL_2 \times
SL_{2m+1}, \mathbb{C}^2 \otimes \Lambda^2 \mathbb{C}^{2m+1})$, $m
\geq 4$, which have infinite number of orbits. In fact the failure
of having a finite number of orbits in these representations is
related to the fact that $PGL_1$ does not act $m$-transitively on
$\mathbb{P}^1$ for $m \geq 4$. In Section \ref{sect:12} we will
discuss their orbit structures in detail.

\medskip
\noindent \textbf{Remark 1}: After the completion of the
preliminary version of this article, we are informed by Landsberg
and later by Rubenthaler that most of our results are already
scattered around in Vinberg \cite{ref23}, Rubenthaler
\cite{ref24}, \cite{ref25} and Kac \cite{ref21}, and are closely
related to Landsberg recent joint works with Manivel and Robles
respectively. More specifically, Theorem \ref{thm:7.1} is similar
to Lemma 1.3 in \cite{ref21} which has been established in Vinberg
\cite{ref23}, and Theorem \ref{thm:6.3} was first proved in
\cite{ref24}. On the other hand, Landsberg and Manivel
\cite{ref22} provide a geometric description via projective
geometry for the minuscule representations, which are subclasses
of representations possessing open orbits. Landsberg also informed
us that his recent joint work with Robles extends the geometric
description in \cite{ref22} to the general case.
\medskip

Let us consider the $GL_n$-case in greater details. There is a
classical result about the irreducible representations of
$GL_n(\mathbb{R})$:
The fundamental representations $\Lambda^k\mathbb{R}^n (k\leq
\frac{n}{2})$ of $GL_n(\mathbb{R})$ has an open orbit if and only
if $(n,k)$ lies in one of the following classes:
\begin{equation} \label{thm:1.1}
\begin{aligned}
\text{(i) }&n \geq 2, k=1;\\
\text{(ii) } &n \geq 4, k=2;\\
\text{(iii) } &n=6,7,8, k=3.
\end{aligned}
\end{equation}

Due to the isomorphism $\Lambda^k\mathbb{R}^n \cong
\Lambda^{n-k}(\mathbb{R}^n)^*$ as $GL_n(\mathbb{R})$
representations, we can confine to the cases where $k \leq
\frac{n}{2}$. One observation is that such configurations can be
reinterpreted as follows: Starting from a Dynkin diagram of type
$A_{n-1}$
\begin{displaymath}
\xymatrix{
\circ\ar@{-}[r]_{\!\!\!\!\!\!\!\!\!\!\!\!\!\!\!\!\!\!\!1} &\circ \ar@{-}[r]_{\!\!\!\!\!\!\!\!\!\!\!\!\!\!\!\!\!\!\!\!2} &\cdots \ar@{-}[r] &\circ \ar@{-}[r]_{\!\!\!\!\!\!\!\!\!\!\!\!\!\!\!\!\!n-2} &\circ \ar@{-}[l]^{\,\,\,\,\,\qquad \qquad n-1}\\
}
\end{displaymath}
one tries to add an extra node to obtain another simply-laced
Dynkin diagram of one higher rank. According to the classification
of reduced root systems \cite{ref8}, we have a full list of
possibilities:
\begin{enumerate}
    \item $A_{n-1} \rightarrow A_n$
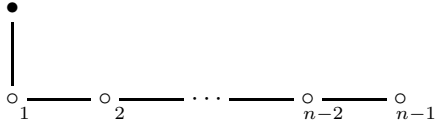
    \begin{displaymath}
\xymatrix{
\bullet \ar@{-}[d]\\
\circ\ar@{-}[r]_{\!\!\!\!\!\!\!\!\!\!\!\!\!\!\!\!\!\!\!1} &\circ \ar@{-}[r]_{\!\!\!\!\!\!\!\!\!\!\!\!\!\!\!\!\!\!\!\!2} &\cdots \ar@{-}[r] &\circ \ar@{-}[r]_{\!\!\!\!\!\!\!\!\!\!\!\!\!\!\!\!\!n-2} &\circ \ar@{-}[l]^{\,\,\,\,\,\qquad \qquad n-1}\\}
\end{displaymath}
    \item $A_{n-1} \rightarrow D_n$
    
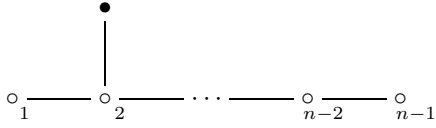
\begin{displaymath}
\xymatrix{
& \bullet \ar@{-}[d] \\
\circ\ar@{-}[r]_{\!\!\!\!\!\!\!\!\!\!\!\!\!\!\!\!\!\!\!1} &\circ \ar@{-}[r]_{\!\!\!\!\!\!\!\!\!\!\!\!\!\!\!\!\!\!\!\!2} &\cdots \ar@{-}[r] &\circ \ar@{-}[r]_{\!\!\!\!\!\!\!\!\!\!\!\!\!\!\!\!\!n-2} &\circ \ar@{-}[l]^{\,\,\,\,\,\qquad \qquad n-1}\\
}
\end{displaymath}
    \item $A_{n-1} \rightarrow E_n$ for $n=6,7,8$
    
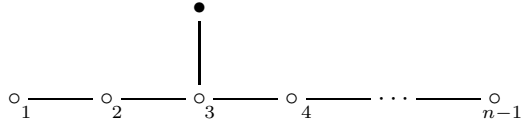
\begin{displaymath}
\xymatrix{
&&\bullet \ar@{-}[d]\\
\circ \ar@{-}[r]_{\!\!\!\!\!\!\!\!\!\!\!\!\!\!\!\!\!\!\!1} &\circ \ar@{-}[r]_{\!\!\!\!\!\!\!\!\!\!\!\!\!\!\!\!\!\!\!\!2} &\circ \ar@{-}[r]_{\!\!\!\!\!\!\!\!\!\!\!\!\!\!\!\!\!\!\!\!3} &\circ \ar@{-}[r]_{\!\!\!\!\!\!\!\!\!\!\!\!\!\!\!\!\!\!\!\!4} &\cdots  &\circ \ar@{-}[l]^{\,\,\,\qquad \qquad n-1}\\
}
\end{displaymath}
\end{enumerate}
In other words, they are obtained from attaching the extra node to
the "$k^{th}$ node" of the original diagram. Then the possible
pairs of $(n,k)$ coincide with the list given in (\ref{thm:1.1}).

The complete dictionary between the existence of open orbits in
$\Lambda^k\mathbb{R}^n$ and the simply-laced extensions of Dynkin
diagram of type $A_{n-1}$ at the "$k^{th}$ node" suggests a
representation-theoretic explanation of this phenomenon. This will
constitute the main content of this paper.

In the case of the fundamental representation
$\Lambda^k(\mathbb{R}^n)^*$ of $GL_n(\mathbb{R})$, an element
$\omega \in \Lambda^k(\mathbb{R}^n)^*$ which lies in an open orbit
is called a stable form. Clearly this notion is independent of the
choice of coordinates and hence it can be defined on any smooth
manifolds. Hitchin \cite{ref19} studied closed differential forms
with such properties. These stable forms have the advantage that
they are stable under deformations and are the critical points of
the associated volume functionals in their respective cohomology
classes, which can be treated as a nonlinear version of the Hodge
theory. The symmetry group $Aut(\mathbb{R}^n,\omega)$ of a stable
form $\omega$ is just the isotropy subgroup of $GL_n(\mathbb{R})$
at $\omega$. For example, when $n$ is even the geometry of stable
two forms is the symplectic geometry and
$Aut(\mathbb{R}^{n},\omega)= Sp(n,\mathbb{R})$. Moreover, the
$E_n$ cases are related to the exceptional geometries, these
geometries are essential to mathematical physics, especially in
developing mathematical models for string theory, these are
studied by Witten \cite{ref20} and his collaborators. For
instance, the geometry of stable $3$-forms on $7$-manifolds are
known as the $G_2$-geometry which is an essential ingredient in
the $M$-theory. In general, given a representation $V$ of a linear
reductive group $G_0$, we try to construct a new semisimple Lie
algebra
$$\mathfrak{g}=\mathfrak{g}_0 \oplus V \oplus V^* \oplus \cdots.$$ But this is
possible only when there is an open orbit, or infinitesimally, we
have $\mathfrak{g}_0= Aut(V,\omega) \oplus V$ for some $\omega \in
V$. It suggests that the representation theory of $\mathfrak{g}$
should be a useful tool to study the geometry of manifolds with
stable forms.


Roughly speaking, we have established the following one-to-one
correspondence:
\begin{align*}
\left\{%
\begin{array}{cc}
    \text{Irreducible prehomogeneous}\\
    \text{vector spaces of parabolic type}\\
\end{array}%
\right\} \longleftrightarrow
\left\{%
\begin{array}{cc}
    \text{Augmentations of}\\
    \text{Dynkin diagrams}\\
\end{array}%
\right\}.
\end{align*}
The term "prehomogeneous vector spaces\,(PVS)" was first
introduced by M. Sato in 1961, since then a lot of results
concerning these objects have been established, in particular the
classification of irreducible prehomogeneous vector spaces was
completed in \cite{ref13}, we will discuss them in Section
\ref{sect:5}. According to a result of Richardson \cite{ref18},
one source of irreducible prehomogeneous vector spaces is obtained
by considering a parabolic subgroup of a complex semisimple Lie
group, from which we obtain a representation of its Levi factor on
the vector space $\mathfrak{u}/ [\mathfrak{u},\mathfrak{u}]$,
where $\mathfrak{u}$ is the nilpotent radical of the corresponding
parabolic subalgebra. Those irreducible prehomogeneous vector
spaces from this origin are said to be of parabolic type. Indeed,
they lies in a subclass of irreducible prehomogeneous vector
spaces which consists of a finite number of orbits. However, they
fail to occupy the whole subclass with a few exceptions which fall
into the class of prehomogeneous vector spaces of twisted affine
type. We will justify our terminology in Section \ref{sect:11}.

Now let's sketch our approach and state our main results. Let $G$
be a connected complex semisimple Lie group with Lie algebra
$\mathfrak{g}$. Upon choosing a Cartan subalgebra $\mathfrak{h}$,
there associates a root system $\Delta$ of $\mathfrak{g}$. Then we
arbitrarily pick up a system of simple roots
$\Pi=\{\alpha_0,\alpha_1,\ldots, \alpha_\ell\}$ and define $c\in
\mathfrak{h}$ to be the unique element such that $\alpha_0(c)=1$
and $\alpha_i(c)=0$ for all $1\leq i \leq \ell$. Set
$\mathfrak{g}_i$ to be the eigenspace of ad $c$ in $\mathfrak{g}$
with eigenvalue $i$, so that we obtain a
$\mathbb{Z}$-gradation\footnote{See Appendix \ref{sect:3} for its
definition.}
$$\mathfrak{g} = \bigoplus_{i \in \mathbb{Z}} \mathfrak{g}_i.$$
It follows immediately that $\mathfrak{g}_0$ is a regular
reductive subalgebra of $\mathfrak{g}$ as the centralizer of $c$
and $\mathfrak{g}$ is called an ambient Lie algebra containing
$\mathfrak{g}_0$. Let $G_0=Z_G(c)^0$ be the closed connected
subgroup of $G$ with Lie algebra $\mathfrak{g}_0$. It is then
clear that all $\mathfrak{g}_i$ are invariant under $G_0$ and thus
are its representations. Note that under the Killing form of
$\mathfrak{g}$, we can identify $\mathfrak{g}_i$ as the dual space
of $\mathfrak{g}_{-i}$ for all $i \neq 0$ and that such
identification is $G_0$-equivariant, i.e. they are dual as
$G_0$-representations. Let $\mathfrak{g}_0^{ss}:= [\mathfrak{g}_0,
\mathfrak{g}_0]$ denote the semisimple part of $\mathfrak{g}_0$.
Then the inclusion $\mathfrak{g}_0^{ss} \subset \mathfrak{g}$
induces a corresponding inclusion of their Dynkin diagrams
$\Gamma(\mathfrak{g}_0^{ss}) \subset \Gamma(\mathfrak{g})$. Our
main result is that to every connected augmentation of Dynkin
diagrams, there exists a unique irreducible reduced PVS with an
extra data called connecting multiplicities to be defined in
Section \ref{sect:6}. More precisely, we have the following
theorem.

\begin{theorem} \label{thm:1.2}
Let $G,G_0, \mathfrak{g}=\bigoplus_{i \in \mathbb{Z}}
\mathfrak{g}_i$ be defined as in the above paragraph, we have:
\begin{enumerate}
    \item For $i \neq 0$, $\mathfrak{g}_i$ are weight multiplicity
free\footnote{See Definition \ref{def:5.1}.} irreducible
representations of $G_0$ with a finite number of orbits; in
particular, it implies that $(G_0,\mathfrak{g}_i)$ are
prehomogeneous vector spaces\footnote{See Definition
\ref{def:5.2}.}.
    \item Every augmentation of Dynkin diagrams can be realized by
    a connected complex semisimple Lie group $G$ with a suitable
    choice of simple root $\alpha_0$, i.e. it can be expressed in
    the form $(\Gamma(\mathfrak{g}),
    \Gamma(\mathfrak{g}_0^{ss}))$.
    \item There is an one-to-one correspodence
    \begin{align*}
\left\{%
\begin{array}{cc}
    \text{Isogeny classes of}\\
    \text{simply-laced irreducible}\\
    \text{prehomogeneous vector spaces}\\
    \text{with finite number of orbits}\\
\end{array}%
\right\} \longleftrightarrow
\left\{%
\begin{array}{cc}
    \text{Simply-laced}\\
    \text{augmentations of}\\
    \text{Dynkin diagrams}\\
\end{array}%
\right\}.
\end{align*}
    Namely, we associate to a simply-laced augmentation of Dynkin diagrams $(\Gamma(\mathfrak{g}),
    \Gamma(\mathfrak{g}_0^{ss}))$ the
    simply-laced prehomogeneous vector space $(G_0,
    \mathfrak{g}_{-1})$.
    \item There is an one-to-one correspondence
    \begin{align*}
\left\{%
\begin{array}{cc}
    \text{Irreducible prehomogeneous}\\
    \text{vector spaces of}\\
    \text{parabolic type together}\\
    \text{with their connecting multiplicities}\\
\end{array}%
\right\} \longleftrightarrow
\left\{%
\begin{array}{cc}
    \text{Connected}\\
    \text{augmentations of}\\
    \text{Dynkin diagrams}\\
\end{array}%
\right\}.
\end{align*}
Explicitly, we assign to each connected augmentation of Dynkin
diagrams $(\Gamma(\mathfrak{g}),
    \Gamma(\mathfrak{g}_0^{ss}))$ the irreducible reduced prehomogeneous
    vector space $$(G_0,
    \mathfrak{g}_{-1}, \nu(\mathfrak{g},\mathfrak{g}_{-1}))$$ with
    the corresponding connecting multiplicities
    $\nu(\mathfrak{g},\mathfrak{g}_{-1})$.
\end{enumerate}
\end{theorem}


Now we illustrate how to apply Theorem \ref{thm:1.2} to our
motivating question at the beginning. The first observation is
that given the representation $\Lambda^k\mathbb{R}^n$ of
$GL_n(\mathbb{R})$, we can complexify it to a representation of
$GL_n(\mathbb{C})$ and then restricted to $SL_n(\mathbb{C})$. By
taking differentials, we obtain the representation
$\Lambda^k\mathbb{C}^n$ of $\mathfrak{sl}_n \mathbb{C}$. Except
the trivial cases where $k=0$ or $n$, for all other cases included
in (\ref{thm:1.1}), $\Lambda^k\mathbb{C}^n$ are corresponding to
the $k^{th}$ or $(n-k)^{th}$ fundamental weight of $\mathfrak{sl}_n
\mathbb{C}$, which is exactly corresponding to the $k^{th}$ or
$(n-k)^{th}$ node of its Dynkin diagram. It provides one possible
linkage between the two sets of objects.

\begin{table}
  \centering
  \caption{Graded pieces of semisimple Lie algebras associated to $(GL_n,\Lambda^k\mathbb{C}^n)$} \label{table:1}
  $\displaystyle{\mathfrak{g}= \bigoplus_{i=-3}^3 \mathfrak{g}_i\,\,,\quad \mathfrak{g}_{-i} \simeq
  \mathfrak{g}_i}$\\
\begin{tabular}{|c|c|c|c|c|}
  \hline
  \phantom{$\displaystyle{\bigoplus}$}$\mathfrak{g}$\phantom{$\displaystyle{\bigoplus}$} & $\mathfrak{g}_0$ & $\mathfrak{g}_{-1}$ & $\mathfrak{g}_{-2}$ & $\mathfrak{g}_{-3} $\\
  \hline
  $\mathfrak{sl}_2 \times \mathfrak{sl}_n$ & $\mathfrak{gl}_n$ & $\mathbb{C}$ & $0$ & \phantom{$\displaystyle{\bigoplus}$}$0$\phantom{$\displaystyle{\bigoplus}$}\\
  $\mathfrak{sl}_{n+1}$ & $\mathfrak{gl}_n$ & $\mathbb{C}^n$ & $0$ & \phantom{$\displaystyle{\bigoplus}$}$0$\phantom{$\displaystyle{\bigoplus}$}\\
  $\mathfrak{so}_{2n}$ & $\mathfrak{gl}_n$ & $\Lambda^2\mathbb{C}^n$ & $0$ & \phantom{$\displaystyle{\bigoplus}$}$0$\phantom{$\displaystyle{\bigoplus}$}\\
  $\mathfrak{e}_6$ & $\mathfrak{gl}_6$ & $\Lambda^3\mathbb{C}^6$ & $\Lambda^6\mathbb{C}^6$ & \phantom{$\displaystyle{\bigoplus}$}$0$\phantom{$\displaystyle{\bigoplus}$}\\
  $\mathfrak{e}_7$ & $\mathfrak{gl}_7$ & $\Lambda^3\mathbb{C}^7$ & $\Lambda^6\mathbb{C}^7$ & \phantom{$\displaystyle{\bigoplus}$}$0$\phantom{$\displaystyle{\bigoplus}$}\\
  \phantom{$\displaystyle{\bigoplus}$}$\mathfrak{e}_8$\phantom{$\displaystyle{\bigoplus}$} & $\mathfrak{gl}_8$ & $\Lambda^3\mathbb{C}^8$ & $\Lambda^6\mathbb{C}^8$ & $\mathbb{C}^8 \otimes \Lambda^8\mathbb{C}^8$\\
  \hline
\end{tabular}
\end{table}

From Table \ref{table:1}, we see that the representation
$\Lambda^k\mathbb{C}^n$ always exists as the $(-1)$-graded
component and it turns out to be the case in general. Then by
Theorem \ref{thm:1.2}, $(GL_n\mathbb{C}, \Lambda^k\mathbb{C}^n)$
have open orbits exactly when $(n,k)$ are as listed in
(\ref{thm:1.1}). Finally by a theorem of Whitney, we successfully
translate the result back to the real cases when the corresponding
complex representation of its real form is of real type; in
particular it is always the cases for split real
forms\footnote{Complete classification of real forms of
irreducible prehomogeneous vector spaces of parabolic type is
obtained in \cite{ref25}}.

Let us briefly describe the content of the paper. In Sections
\ref{sect:5} and \ref{sect:6}, we will set up the general
framework and the terminology used throughout this paper. The
main result on the finiteness of orbits will be established in Section
\ref{sect:7}. Then the termination of $\mathbb{Z}$-gradations will
be discussed in Section \ref{sect:8}. After that, we will give an
explicit construction of generic elements in the simply-laced
cases in Section \ref{sect:9}. The proof of Theorem \ref{thm:1.2}
will be completed in Section \ref{sect:10}. Sections \ref{sect:11}
and \ref{sect:12} are devoted to the discussion of the two
exceptional cases in our construction. Finally the first two
appendices present the basics of $\mathbb{Z}$-gradations and
algebraic groups, and the tables are contained in the last
appendix.

\section{Weight Multiplicity Free Representations and Prehomogeneous Vector Spaces} \label{sect:5}
In this section, we will introduce two notions in representation
theory which have been well understood for a long time, namely
that of weight multiplicity free representations and prehomogeneous
vector spaces. Both objects have been completely classified and
proved to be useful in many branches of mathematics. Here we will
use them to give a necessary condition for the existence of
augmentation of Dynkin diagrams.

\begin{definition} \label{def:5.1}
Let $\mathfrak{g}$ be a complex semisimple Lie algebra. A representation $V$ of $\mathfrak{g}$ is said to be weight multiplicity free if every weight space is one dimensional.
\end{definition}

The classification of irreducible weight multiplicity free representations of complex simple Lie algebras can be found in \cite{ref5}. The complete list is as follows:
\begin{enumerate}
    \item $A_\ell\,(\mathfrak{sl}_{\ell +1}\mathbb{C})$:
\begin{enumerate}
    \item The fundamental representations $\Lambda^m\mathbb{C}^{\ell+1}$ with highest weight $\omega_m$, for $m=1, \ldots, \ell$.
    \item The symmetric tensor powers $S^m\mathbb{C}^{\ell +1}$ and $S^m(\mathbb{C}^{\ell +1})^*$ with highest weights $m\omega_1$ and $m\omega_\ell$, for $m \in \mathbb{Z}_{\geq 0}$.
\end{enumerate}
    \item $B_\ell\,(\mathfrak{so}_{2\ell +1}\mathbb{C})$:
\begin{enumerate}
    \item The standard representation $\mathbb{C}^{2\ell+1}$ with highest weight $\omega_1$.
    \item The spin representation $S$ with highest weight $\omega_\ell$.
\end{enumerate}
    \item $C_\ell\,(\mathfrak{sp}_{2\ell}\mathbb{C})$:
\begin{enumerate}
    \item The standard representation $\mathbb{C}^{2\ell}$ with highest weight $\omega_1$.
    \item When $\ell=2$ or $3$, the last fundamental representation, $\Lambda^2_{\text{prim}}\mathbb{C}^4$ and $\Lambda^3_{\text{prim}}\mathbb{C}^6$ respectively, with highest weight $\omega_\ell$.
\end{enumerate}
    \item $D_\ell\,(\mathfrak{so}_{2\ell}\mathbb{C})$:
\begin{enumerate}
    \item The standard representation $\mathbb{C}^{2\ell}$ with highest weight $\omega_1$.
    \item The two half-spin representations $S^+$ and $S^-$ with highest weights $\omega_{\ell-1}$ and $\omega_\ell$ respectively.
\end{enumerate}
    \item $E_\ell\,(\ell=6,7,8)$:
\begin{enumerate}
    \item The two $27$-dimensional representations of $E_6$ with highest weights $\omega_1$ and $\omega_6$.
    \item The $56$-dimensional representations of $E_7$ with highest weight $\omega_7$.
    \item There are no weight multiplicity free representations for $E_8$.
\end{enumerate}
    \item $F_4$:
    There are no weight multiplicity free representations for $F_4$.
    \item $G_2$:
    The $7$-dimensional representation of $G_2$ with highest weight $\omega_2$.
\end{enumerate}
Here the numbering of the fundamental weights $\omega_i$ are adopted to that of Bourbaki \cite{ref2}(also see Table \ref{table:2}).

\begin{proposition}
Let $\mathfrak{g}= \mathfrak{a}_1 \times  \cdots \times \mathfrak{a}_k$ be a decomposition of a semisimple Lie algebra $\mathfrak{g}$ into simple ideals $\mathfrak{a}_i, i = 1, \ldots ,k$, and suppose that $V_i$ is a finite dimensional representation of $\mathfrak{a}_i$ for each $i$. Then $V_1 \otimes \cdots \otimes V_k$ is a weight multiplicity free representation of $\mathfrak{g}$ if and only if all the $V_i$'s are weight multiplicity free.
\end{proposition}

\begin{proposition} \label{prop:5.1}
Let $V$ be a weight multiplicity free representation of a complex semisimple Lie algebra $\mathfrak{g}$. If all weights of $V$ are congruent modulo the root lattice $\Lambda$ of $\mathfrak{g}$, then $V$ is irreducible.
\end{proposition}
\begin{proof}
Fix a Cartan subalgebra $\mathfrak{h}$ and a system of simple roots $\Pi=\{\alpha_1, \ldots, \alpha_\ell\}$ so that $\Lambda=\mathbb{Z}\left\langle \alpha_1, \ldots, \alpha_\ell\right\rangle$, and let $H_i$ be the unique element in $[\mathfrak{g}_{\alpha_i},\mathfrak{g}_{-\alpha_i}]$ with $\alpha_i(H_i)=2$. It suffices to show that the highest weight of $V$ is unique. Suppose on the contrary that $\lambda$ and $\mu$ are two distinct highest weights of $V$. Since $\lambda$ and $\mu$ are congruent modulo $\Lambda$, $\lambda-\mu=\alpha$ for some $\alpha \in \Lambda$. Now by separating the positive and negative parts of $\alpha$ as a linear combination of $\alpha_1, \ldots, \alpha_\ell$, we can find two disjoint subset $I,J$ of $\{\alpha_1, \ldots, \alpha_\ell\}$ and positive integers $n_i,m_j$ for every $i\in I, j \in J$, such that $$\eta:= \lambda- \sum_{i \in I}n_i \alpha_i = \mu - \sum_{j \in J}m_j \alpha_j.$$

Note that for each $i \notin J$, $$\eta(H_i)= \mu(H_i)-\sum_{j \in J} m_j \alpha_j(H_i) \geq 0$$ since $\alpha_j(H_i)\leq 0$ for all $j \neq i$. Similarly, we have for each $j \in I$, $$\eta(H_j)=\lambda(H_j)- \sum_{i\in I}n_i \alpha_i(H_j) \geq 0.$$ As $I$ and $J$ are disjoint, we conclude that $\eta(H_i)\geq 0$ for all $i = 1,\ldots, \ell$ which implies that $\eta$ is a dominant integral weight in $\mathfrak{h}^*$. It follows that $\eta$ must be a weight of $V$ with multiplicity at least two since each of the highest weight submodules of $V$ of weights $\lambda$ and $\mu$ contributes at least one dimension to the weight space $V_\eta$. But this contradicts that $V$ is weight multiplicity free.
\end{proof}

\begin{definition} \label{def:5.2}
Let $G$ be a complex reductive algebraic group and $V$ a rational representation of $G$. $(G,V)$ is called a prehomogeneous vector space if there exists a dense $G$-orbit in $V$.
\end{definition}

According to Proposition \ref{prop:4.1} in Appendix \ref{sect:4},
the dense orbit must be open. In other words, prehomogeneous
vector spaces are just representations with exactly one open
orbit. All irreducible prehomogeneous vector spaces have been classified
in \cite{ref13}.

\begin{proposition} \label{prop:5.2}
Given a representation $V$ of a linear connected algebraic group $G$. Then the following conditions are equivalent:
\begin{enumerate}
    \item $(G,V)$ is a prehomogeneous vector space.
    \item There exists a vector $v \in V$ such that $\dim G_v = \dim G - \dim V$, where $G_v=\{g \in G|g\cdot v=v\}$.
    \item There exists a vector $v \in V$ such that $\mathfrak{g}\cdot v = V$.
\end{enumerate}
\end{proposition}

\begin{theorem}[Richardson \cite{ref18}] \label{thm:5.3}
Let $G$ be a connected complex semisimple Lie group, and let $P$
be a parabolic subgroup of $G$ with Levi decomposition $P=LU$,
where $L$ is its Levi factor and $U$ its unipotent radical. If
$\mathfrak{u}=Lie(U)$, then $(L, \mathfrak{u}/
[\mathfrak{u},\mathfrak{u}])$ is a prehomogeneous vector space.
\end{theorem}

\begin{definition}
Let $(G_0,V)$ be a prehomogeneous vector space, where $G_0$ is a
connected complex reductive Lie group of the form $GL_1 \times
G_0^{ss}$ for some complex semisimple Lie group $G_0^{ss}$.
\begin{enumerate}
    \item $(G_0,V)$ is said to of parabolic
type if it can be obtained from a connected complex semisimple Lie
group $G$ in the sense of Theorem \ref{thm:5.3}.
    \item $(G_0,V)$ is said to be reduced if its dimension is
    minimal over all prehomogeneous vector spaces which are
    castling equivalent\footnote{$(GL_n \times G, \mathbb{C}^n \otimes V) \cong (GL_{n-m} \times G, \mathbb{C}^{n-m} \otimes V^*)$ where $\dim V = m < n$ induces an equivalence relation on the set of PVS, which is called the castling equivalence.} to it.
    \item $(G_0,V)$ is said to be irreducible if $V$ is
    irreducible as a $G_0$-representation.
\end{enumerate}
\end{definition}

\section{Augmentations of Dynkin Diagrams} \label{sect:6}
It is well-known that every complex semisimple Lie algebra $\mathfrak{g}$ admits a unique Dynkin diagram $\Gamma(\mathfrak{g})$ determined by the associated Cartan matrix, up to permutations of numbering of its entries. We also know that $\Gamma(\mathfrak{g})$ is connected if and only if $\mathfrak{g}$ is simple. In this section, we would like to study when one can add an extra node to a given Dynkin diagram with corresponding relations of the edges attaching to the node so that it remains a Dynkin diagram of some semisimple Lie algebra of higher rank. In other words, we want to study all pairs of Dynkin diagrams $(\Gamma,\Gamma_0)$ consisting of a Dynkin diagram $\Gamma$ and a subdiagram $\Gamma_0$ obtained by removing a single node and all edges attached to it.

\begin{definition}
An augmentation of Dynkin diagrams is a pair $(\Gamma,\Gamma_0)$ of Dynkin diagrams such that $\Gamma_0$ is a subdiagram of $\Gamma$ obtained by removing exactly one node and all the edges connected to it.
\end{definition}

Note that subdiagrams of a Dynkin diagram correspond exactly to
the principal minors of the corresponding Cartan matrix. Therefore
any subdiagram of a Dynkin diagram is also a Dynkin diagram and
the above definition makes sense. To represent an augmentation of
Dynkin diagrams $(\Gamma,\Gamma_0)$ diagrammatically, we will use
the Dynkin diagram $\Gamma$ with a painted node indicating the
omitted node in $\Gamma_0$.

Starting with a semisimple Lie algebra $\mathfrak{g}$, our approach is to give a realization of $\Gamma_0$ as a subsystem of simple roots of a semisimple subalgebra of $\mathfrak{g}$ through a $\mathbb{Z}$-gradation of $\mathfrak{g}$, and we will associate to it a collection of irreducible representations which detect the validity of such pair.

First of all, let's set up some notations. Let $G$ be a connected
complex semisimple Lie group with Lie algebra $\mathfrak{g}$, and
$\mathfrak{h}$ be a Cartan subalgebra of $\mathfrak{g}$. Then we
have a root space decomposition of $\mathfrak{g}$ with respect to
$\mathfrak{h}$ $$\mathfrak{g}= \mathfrak{h} \oplus
\bigoplus_{\alpha \in \Delta} \mathfrak{g}_\alpha$$ where $\Delta$
is the root system of $\mathfrak{g}$ with respect to
$\mathfrak{h}$. Assume that rank$(\mathfrak{g})= \ell + 1$, and
let $\Pi=\{\alpha_0,\alpha_1, \ldots, \alpha_\ell\}\subset \Delta$
be a system of simple roots. For each $i \in \mathbb{Z}$, set
$\Delta_i= (\mathbb{Z}\left\langle \alpha_1, \ldots, \alpha_\ell
\right\rangle + i \alpha_0)\cap \Delta$ so that $$\Delta=
\bigsqcup_{i \in \mathbb{Z}} \Delta_i.$$ Now for each $i \neq 0$,
denote $\displaystyle{\mathfrak{g}_i= \bigoplus_{\alpha \in
\Delta_i} \mathfrak{g}_\alpha}$, and define
$\displaystyle{\mathfrak{g}_0= \mathfrak{h} \oplus
\bigoplus_{\alpha \in \Delta_0} \mathfrak{g}_\alpha}$. Then it is
easy to verify that $\displaystyle{\mathfrak{g}= \bigoplus_{i \in
\mathbb{Z}}\mathfrak{g}_i}$ is a $\mathbb{Z}$-gradation. From
Proposition \ref{prop:3.1} in Appendix \ref{sect:3},
$\mathfrak{g}_0$ is a reductive subalgebra of $\mathfrak{g}$ and
thus $\mathfrak{g}_0^{ss}$ is a semisimple subalgebra of rank
$\ell$.

Let $c$ be the element in $\mathfrak{h}$ such that $\alpha_0(c)=1$ and $\alpha_i(c)=0$ for all $1 \leq i \leq \ell$. Then we can write $\mathfrak{h}= \mathbb{C}c \oplus \mathfrak{t}$ where $\mathfrak{t}$ is the orthogonal complement of $\mathbb{C}c$ in $\mathfrak{h}$ with respect to the Killing form of $\mathfrak{g}$. It is clear that $\mathfrak{t}$ is a Cartan subalgebra of $\mathfrak{g}_0^{ss}$ and the corresponding root system $\widetilde{\Delta_0}= \{\alpha|_\mathfrak{t}\, | \alpha \in \Delta_0\}$ has a subsystem of simple roots given by $\widetilde{\Pi_0}=\{\alpha_i|_\mathfrak{t}\, | 1 \leq i \leq \ell \}$.

If we identify $\widetilde{\Pi_0}$ with $\Pi_0=\{\alpha_1, \ldots, \alpha_\ell\}$, we see that $(\Gamma(\mathfrak{g}), \Gamma(\mathfrak{g}_0^{ss}))$ is an augmentation of Dynkin diagrams. As long as only augmentations of Dynkin diagrams are concerned, the choices of the Cartan subalgebra $\mathfrak{h}$ and the system of simple roots $\Pi$ are inessential, and we will fix $\mathfrak{h}$ and $\Pi$ once and for all in the remaining part of this paper.

The element $c$ constructed above plays an important role in the
structure of the $\mathbb{Z}$-gradation of $\mathfrak{g}$, for
instance, we have $\mathfrak{g}_i=\{X \in \mathfrak{g}|
[c,X]=iX\}$ for all $i \in \mathbb{Z}$; in particular,
$\mathfrak{g}_0$ is the centralizer of $c$ in $\mathfrak{g}$ with
center $\mathfrak{z}_{\mathfrak{g}_0}=\mathbb{C}c$.

\begin{lemma} \label{lemma:6.1}
Let $\displaystyle{\mathfrak{g}= \bigoplus_{i \in \mathbb{Z}} \mathfrak{g}_i}$ be the $\mathbb{Z}$-gradation constructed as in the above discussion, and let $\left\langle \cdot,\cdot\right\rangle_\mathfrak{h}$ and $\left\langle \cdot , \cdot \right\rangle_\mathfrak{t}$ be the Cartan products on $\mathfrak{h}^*$ and $\mathfrak{t}^*$ respectively.
Then $\left\langle \alpha,\beta\right\rangle_\mathfrak{h} = \left\langle \alpha|_\mathfrak{t}, \beta|_\mathfrak{t} \right\rangle_\mathfrak{t}$ for all $\alpha \in \Delta, \beta \in \Delta_0$.
\end{lemma}
\begin{proof}
Let $H_\beta$ be the unique element in $[\mathfrak{g}_\beta,\mathfrak{g}_{-\beta}] \subset \mathfrak{t}$ such that $\beta(H_\beta)=2$. Since $\mathfrak{g}_\beta= (\mathfrak{g}_0^{ss})_{\beta|_\mathfrak{t}}$, $H_\beta$ is also the unique element in $[(\mathfrak{g}_0^{ss})_{\beta|_\mathfrak{t}},(\mathfrak{g}_0^{ss})_{-\beta|_\mathfrak{t}}]$ such that $\beta|_\mathfrak{t}(H_\beta)=2$. It follows that $$\left\langle \alpha|_\mathfrak{t}, \beta|_\mathfrak{t} \right\rangle_\mathfrak{t}= \alpha|_\mathfrak{t}(H_\beta) = \alpha(H_\beta) = \left\langle \alpha,\beta\right\rangle_\mathfrak{h}.$$
\end{proof}

\begin{theorem} \footnote{The theorem was first proved in \cite{ref24}. See also Remark 1.}  \label{thm:6.3}
Let $\displaystyle{\mathfrak{g}= \bigoplus_{i \in \mathbb{Z}}
\mathfrak{g}_i}$ be defined as above.
\begin{enumerate}
    \item For each $k \neq 0$, $\mathfrak{g}_k$ is an irreducible weight multiplicity free representation of $\mathfrak{g}_0^{ss}$.
    \item $-\alpha_0|_\mathfrak{t}$ is the highest weight of $\mathfrak{g}_{-1}$ as a $\mathfrak{g}_0^{ss}$-representation.
\end{enumerate}
\end{theorem}
\begin{proof}
\begin{enumerate}
    \item Note that from $\displaystyle{\mathfrak{g}_k= \bigoplus_{\alpha \in \Delta_k} \mathfrak{g}_\alpha}$, we have each root space $\mathfrak{g}_\alpha$ in $\mathfrak{g}_k$ is contained in the weight space $\mathfrak{g}_k$ as a $\mathfrak{g}_0^{ss}$-module
of weight $\alpha|_\mathfrak{t}$. Now for any two distinct $\alpha, \beta \in \Delta_k$, we have $\alpha- \beta \in \mathbb{Z}\left\langle \alpha_1, \ldots , \alpha_\ell\right\rangle$ and thus $\alpha|_\mathfrak{t}-\beta|_\mathfrak{t}$ is a nonzero element in the root lattice $\Lambda_0$ generated by $\widetilde{\Delta_0}$. Hence distinct root spaces in $\mathfrak{g}_k$ lie in different weight spaces of $[\mathfrak{g}_0, \mathfrak{g}_0]$; in other words, $\displaystyle{\mathfrak{g}_k= \bigoplus_{\alpha \in \Delta_k} \mathfrak{g}_\alpha}$ is precisely the weight space decomposition as $\mathfrak{g}_0^{ss}$-module. Now since each root space is one dimensional, we conclude that $\mathfrak{g}_k$ is a weight multiplicity free representation of $\mathfrak{g}_0^{ss}$. Finally as all weights are congruent to each other modulo the root lattice $\Lambda_0$, $\mathfrak{g}_k$ is irreducible according to Proposition \ref{prop:5.1}.
\item Note that from the proof of (i), we have the set of weights of $\mathfrak{g}_{-1}$ being the restriction of the elements in $\Delta_{-1}$ to $\mathfrak{t}$. Now for each $\alpha \in \Delta_{-1}$, $$\alpha= -\alpha_0 - \sum_{i=1}^\ell n_i \alpha_i$$ for some non-negative integers $n_i(i=1, \ldots , \ell)$. It follows that $$-\alpha_0- \alpha = \sum_{i=1}^\ell n_i \alpha_i$$ is positive in the lexicographical ordering. Thus $-\alpha_0|_\mathfrak{t}$ is the highest weight of $\mathfrak{g}_{-1}$ as a $\mathfrak{g}_0^{ss}$-representation.
\end{enumerate}
\end{proof}

Up to now, we have established (i) and (ii) of Theorem
\ref{thm:1.2}. From Theorem \ref{thm:6.3}, we observe that the
$\mathfrak{g}_0^{ss}$-representation $\mathfrak{g}_{-1}$ imposes
severe constraint on the possible augmentations of Dynkin diagrams
as it gives a finite list of possible weights
$-\alpha_0|_\mathfrak{t}$. Assuming the existence of such
augmentation of Dynkin diagrams, in virtue of Lemma
\ref{lemma:6.1}, $\mathfrak{g}_{-1}$ determines the Cartan matrix
of the possible $\mathfrak{g}$ up to the choice of the values
$\left\langle \alpha_i,\alpha_0\right\rangle_\mathfrak{h}$ for $1
\leq i \leq \ell$. But according to the properties of Cartan
matrices, the only ambiguity happens for those $i$ where
$\left\langle \alpha_0,\alpha_i\right\rangle_\mathfrak{h} =-1$.
Therefore, if we restrict to only simply-laced simple Lie algebra
$\mathfrak{g}$, $\mathfrak{g}_{-1}$ determines completely the
structure of $\mathfrak{g}$.

An alternative method to remove the ambiguity is to associate the
missing vector of integers $(-\left\langle
\alpha_1,\alpha_0\right\rangle_\mathfrak{h}, \ldots, -\left\langle
\alpha_\ell,\alpha_0\right\rangle_\mathfrak{h})$ to the
representation $\mathfrak{g}_{-1}$. Since every Dynkin diagram
does not have any cycles, the extra node can only connected to
each component of $\Gamma(\mathfrak{g}_0^{ss})$ at no more than
one node $\alpha_i$, and only those $\left\langle
\alpha_i,\alpha_0\right\rangle_\mathfrak{h}$ give us information.
Hence we define the notion of connecting multiplicities to capture
this piece of information.

\begin{definition}
Let $\displaystyle{\mathfrak{g}= \bigoplus_{i \in \mathbb{Z}}
\mathfrak{g}_i}$ be defined as above. Suppose that
$\mathfrak{g}_0^{ss}= \mathfrak{a}_1 \times \cdots \times
\mathfrak{a}_r$ is the decomposition of $\mathfrak{g}_0^{ss}$ into
simple ideals. The $\mathbb{Z}$-valued vector
$$\nu(\mathfrak{g},\mathfrak{g}_{-1})=(a_1, \ldots, a_r)
\text{ with } \displaystyle{a_i= \max_{\alpha \in
\Gamma(\mathfrak{a}_i)} -\left\langle
\alpha,\alpha_0\right\rangle_\mathfrak{h}},$$ is called the
connecting multiplicities of $\mathfrak{g}_0^{ss}$-representation
$\mathfrak{g}_{-1}$ with respect to $\mathfrak{g}$.
\end{definition}

\noindent \textbf{Remark 2}: To avoid ambiguity in defining
$\nu(\mathfrak{g},\mathfrak{g}_{-1})$, we adopt the convention
that the simple ideals $\mathfrak{a}_i$ are lined up in
alphabetical order according to their Lie types and among those
with the same Lie type we write the one with smaller rank in
front. If it happens that some of them are exactly the same, we
simply put the values $a_i$ in descending order.
\medskip

Clearly, elements in $\nu(\mathfrak{g},\mathfrak{g}_{-1})$ takes values only
from $\{0,1,2,3\}$ and that $2,3$ cannot appear twice or at the
same time. Also $\mathfrak{g}$ is simple if and only if
$\nu(\mathfrak{g},\mathfrak{g}_{-1})$ does not contain $0$.

The following proposition captures some important direct consequences from the Cartan matrix of $\mathfrak{g}$.

\begin{proposition} \label{prop:6.4}
Let $\displaystyle{\mathfrak{g}= \bigoplus_{i \in \mathbb{Z}}
\mathfrak{g}_i}$ be defined as above, $\omega$ be the highest
weight of $\mathfrak{g}_{-1}$ as a $\mathfrak{g}_0^{ss}$-representation and $A$ be the Cartan matrix of $\mathfrak{g}$
with respect to the system of simple roots $\Pi=\{\alpha_0, \ldots
\alpha_\ell\}$. Suppose $\nu(\mathfrak{g},\mathfrak{g}_{-1})=(a_1,
\ldots, a_r)$ and $a_i = \left\langle
\alpha_i,\alpha_0\right\rangle_\mathfrak{h}$.
\begin{enumerate}
    \item $a_i=0$ if and only if $\left\langle
    \omega,\alpha_i|_{\mathfrak{t}}\right\rangle_\mathfrak{t}=0$.
    \item For all $a_i \neq 0$, we have
    $$\frac{a_i(H_{\alpha_i},H_{\alpha_i})}{\omega(H_{\alpha_i})}$$
    is a nonzero constant independent of $i$.
    \item All principal minors of $A$ are positive definite.
\end{enumerate}
\end{proposition}

In fact, we will show in Section \ref{sect:10} that those conditions in Proposition \ref{prop:6.4} are the only conditions required to construct back the ambient Lie algebra $\mathfrak{g}$. As a result, we abstractly define the connecting multiplicities of an arbitrary irreducible representation.

\begin{definition}
Let $V$ be an irreducible representation of a semisimple Lie
algebra $\mathfrak{g}$ with highest weight $\omega$. A
$\mathbb{Z}$-valued vector $\nu=(a_1, \ldots, a_r)$ is called the
connecting multiplicities of $V$ if the following conditions are
satisfied: There is a system of simple roots $\Pi=\{\alpha_1,
\ldots ,
    \alpha_\ell\}$ so that
\begin{enumerate}
    \item
    $a_i=0$ if and only if $\left\langle
    \omega,\alpha_i\right\rangle=0$.
    \item For all $a_i \neq 0$, we have
    $$\frac{a_i(H_{\alpha_i},H_{\alpha_i})}{\omega(H_{\alpha_i})}$$
    is a nonzero constant independent of $i$.
    \item All principal minors of
    $$A=\left(%
\begin{array}{ccccccc}
  2 & -\langle\omega, \alpha_1\rangle & \cdots & -\langle\omega, \alpha_r\rangle & 0 & \cdots & 0 \\
  -a_1 & \langle\alpha_1, \alpha_1\rangle & \cdots & \langle\alpha_1, \alpha_r\rangle & \langle\alpha_1, \alpha_{r+1}\rangle & \cdots & \langle\alpha_1, \alpha_\ell\rangle \\
  \vdots &  &  &  &  &  &  \\
  -a_r & \vdots &  &  &  &  & \vdots \\
  0 & \vdots &  &  &  &  & \vdots \\
  \vdots &  &  &  &  &  &  \\
  0 & \langle\alpha_\ell, \alpha_1\rangle &  & \cdots & \cdots &  & \langle\alpha_\ell, \alpha_\ell\rangle \\
\end{array}%
\right).$$
    are positive definite.
\end{enumerate}
\end{definition}

\noindent \textbf{Remark 3}: From the definition of $\nu$, $A$ is
a Cartan matrix.
\medskip

There is one further possible reduction of the problem, namely it
suffices to consider simple Lie algebra $\mathfrak{g}$. The reason
is that only the connected component of $\Gamma(\mathfrak{g})$
containing $\alpha_0$ is sensitive to our $\mathbb{Z}$-gradation
of $\mathfrak{g}$. Note that the root system of $\mathfrak{g}$
decomposes into irreducible subsystems which are mutually
orthogonal to each other, each of which corresponds to a connected
component of $\Gamma(\mathfrak{g})$. Thus the root spaces
$\mathfrak{g}_\alpha$ are contained in $\mathfrak{g}_0$ for those
$\alpha$ lying in an irreducible subsystem not containing
$\alpha_0$, and they act trivially on each $\mathfrak{g}_k (k \neq
0)$. Therefore, if $\Gamma_0(\mathfrak{g})$ denotes the connected
component of $\Gamma(\mathfrak{g})$ containing $\alpha_0$ and
$\Gamma_0(\mathfrak{g}_0^{ss})$ is the subdiagram of
$\Gamma(\mathfrak{g}_0^{ss})$ obtained by deleting the nodes lying
in the connected components not containing $\alpha_0$, then
$(\Gamma_0(\mathfrak{g}), \Gamma_0(\mathfrak{g}_0^{ss}))$ is still
an augmentation of Dynkin diagrams, and all the nonzero graded
pieces $\mathfrak{g}_k$ remain unchange.

In the next section, we will show that such $\mathfrak{g}_k$ are
prehomogeneous vector spaces with respect to a closed connected
reductive algebraic group $G_0$ corresponding to the Lie algebra
$\mathfrak{g}_0$.

\section{Orbit Finiteness and Prehomogeneity} \label{sect:7}
Up to now, we have shown that $\mathfrak{g}_k$ for $k \neq 0$ are irreducible weight multiplicity free representations of $\mathfrak{g}_0^{ss}$. By taking into account of the action of the closed connected subgroup $G_0$ of $G$ with Lie algebra $\mathfrak{g}_0$, we succeed in showing that $(G_0,\mathfrak{g}_k)$ are prehomogeneous vector spaces for all $k \neq 0$. Essentially the proof will be separated into two steps: 1) To establish an orbit finiteness statement of $G_0$ on $\mathfrak{g}_k$; 2) to show that there is exactly one open orbit in $\mathfrak{g}_k$, which is the restriction of a nilpotent orbit in $\mathfrak{g}$ onto $\mathfrak{g}_k$.

\begin{theorem} \footnote{The theorem was first proved in \cite{ref23}. See also Remark 1.} \label{thm:7.1}
Let $G$ be a connected complex semisimple Lie group with Lie
algebra $\mathfrak{g}$. Suppose that $\mathfrak{g}$ has a
$\mathbb{Z}$-gradation $\mathfrak{g}= \bigoplus_{i \in \mathbb{Z}}
\mathfrak{g}_i$. Let $G_0$ be the closed connected subgroup of $G$
with Lie algebra $\mathfrak{g}_0$. Then for each $k \neq 0$, the
action of $G_0$ on $\mathfrak{g}_k$ has a finite number of orbits.
\end{theorem}

Before going into the proof, we need a simple lemma.

\begin{lemma} \label{lemma:7.2}
Under the same conditions as in Theorem \ref{thm:7.1}, for each $k\neq 0$, every element in $\mathfrak{g}_k$ is nilpotent in $\mathfrak{g}$.
\end{lemma}
\begin{proof}
Pick any $X \in \mathfrak{g}_k$, we have, for all $i \in
\mathbb{Z}$, ad$_X^m(\mathfrak{g}_i)\subset \mathfrak{g}_{i+mk}$.
Since $\mathfrak{g}$ is finite dimensional and $k \neq 0$,
$\mathfrak{g}_{i+mk}=0$ for sufficiently large $m$. Hence
ad$_X^m$=0 and $X$ is nilpotent in $\mathfrak{g}$.
\end{proof}

\begin{proof}[Proof of Theorem \ref{thm:7.1}.\!\!]
First note that $[\mathfrak{g}_0,\mathfrak{g}_k] \subset
\mathfrak{g}_k$, and hence $\mathfrak{g}_k$ is $G_0$-invariant.
Let $\widetilde{\mathfrak{g}_k}:=$ Ad$(G) \cdot \mathfrak{g}_k$ be
the $G$-saturation of $\mathfrak{g}_k$ in $\mathfrak{g}$. By Lemma
\ref{lemma:7.2}, every element in $\mathfrak{g}_k$ is nilpotent in
$\mathfrak{g}$, so that the same is true for
$\widetilde{\mathfrak{g}_k}$. In other words,
$\widetilde{\mathfrak{g}_k}$ is a union of nilpotent $G$-orbits in
$\mathfrak{g}$, which must be finite since there are only finite
number of nilpotent $G$-orbits in $\mathfrak{g}$ \cite{ref3}. It
remains to show that for every $G$-orbit $\mathcal{O}$ in
$\widetilde{\mathfrak{g}_k}$, $\mathcal{O}\cap \mathfrak{g}_k$
splits into a finite number of $G_0$-orbits in $\mathfrak{g}_k$.

By the definition of $\widetilde{\mathfrak{g}_k}$, there exists $X
\in \mathfrak{g}_k$ such that $\mathcal{O}=$ Ad$(G)\cdot X$, so
that $\mathcal{O}\cap \mathfrak{g}_k \neq \varnothing$. Now for
any $X' \in \mathcal{O}\cap \mathfrak{g}_k$,
$$T_{X'}(\mathcal{O}\cap \mathfrak{g}_k) \subset
[\mathfrak{g},X']\cap
\mathfrak{g}_k=[\mathfrak{g}_0,X']=T_{X'}(\text{Ad}(G_0)\cdot
X').$$ But on the other hand, we have $\mathcal{O}\cap
\mathfrak{g}_k \supset$ Ad$(G_0)\cdot X'$ since $\mathfrak{g}_k$
is $G_0$-invariant. It follows that $T_{X'}(\mathcal{O}\cap
\mathfrak{g}_k)=T_{X'}(\text{Ad}(G_0)\cdot X')$, thus $X'$ is a
nonsingular point of $\mathcal{O}\cap \mathfrak{g}_k$ and
$\text{Ad}(G_0)\cdot X'$ is open in $\mathcal{O}\cap
\mathfrak{g}_k$. As $X' \in \mathcal{O}\cap \mathfrak{g}_k$ is
arbitrary, $\text{Ad}(G_0)\cdot X'$ is also closed in
$\mathcal{O}\cap \mathfrak{g}_k$ for its complement is a union of
such orbits. Thus the $G_0$-orbits in $\mathcal{O}\cap
\mathfrak{g}_k$ are precisely all the connected components of
$\mathcal{O}\cap \mathfrak{g}_k$ and so $\mathcal{O}\cap
\mathfrak{g}_k$ being a smooth manifold can possess only finite
number of $G_0$-orbits.
\end{proof}

\noindent \textbf{Remark 4}: The analogous statement of Theorem
\ref{thm:7.1} for $\mathbb{Z}_m$-gradation
$\displaystyle{\mathfrak{g}= \bigoplus_{i \in \mathbb{Z}_m}
\mathfrak{g}_i}$ holds true as long as $\mathfrak{g}_k,k \neq 0,$
are contained in the nilpotent cone of $\mathfrak{g}$. The line of
proof runs exactly the same except one must replace Lemma
\ref{lemma:7.2} by the above assumption.

\begin{theorem}
Under same conditions as in Theorem \ref{thm:7.1},
$\mathfrak{g}_k$ has a unique open $G_0$-orbit of the form
$\mathcal{O}_k \cap \mathfrak{g}_k$, where $\mathcal{O}_k$ is the
unique open $G$-orbit  in $\widetilde{\mathfrak{g}_k}$.
In particular, $(G_0,\mathfrak{g}_k)$ is a prehomogeneous vector
space for every $k \neq 0$.
\end{theorem}
\begin{proof}
Note that $\mathfrak{g}_k$ is irreducible as an affine variety, which forces all open $G_0$-orbits in $\mathfrak{g}_k$ to be dense and thus coincide. Therefore, $(G_0,\mathfrak{g}_k)$ is trivially a prehomogeneous vector space in virtue of Theorem \ref{thm:7.1}. This proves the second statement. To establish the first assertion, we need the following two lemmas:

\begin{lemma} \label{lemma:7.4}
With the same notations as in Theorem \ref{thm:7.1}, there is a unique nilpotent $G$-orbit in $\widetilde{\mathfrak{g}_k}=$ Ad$(G)\cdot \mathfrak{g}_k$ which is open in $\widetilde{\mathfrak{g}_k}$ for every $k \neq 0$.
\end{lemma}
\begin{proof}
Suppose on the contrary that there are two such nilpotent $G$-orbits $\mathcal{O}', \mathcal{O}''$. Then $\mathcal{O}' \cap \mathfrak{g}_k$ and $\mathcal{O}'' \cap \mathfrak{g}_k$ are nonempty and open in $\mathfrak{g}_k$. As $\mathfrak{g}_k$ is an affine space, and hence irreducible, $\mathcal{O}'\cap \mathfrak{g}_k$ and $\mathcal{O}'' \cap \mathfrak{g}_k$ are dense in $\mathfrak{g}_k$. It follows that $\mathcal{O}'$ intersects with $\mathcal{O}''$ nontrivially, which forces $\mathcal{O}'=\mathcal{O}''$.
\end{proof}

\begin{lemma} \label{lemma:7.5}
Let $\mathcal{O}_k$ be the unique nilpotent $G$-orbit contained in $\widetilde{\mathfrak{g}_k}$ obtained in Lemma \ref{lemma:7.4}. Then $\mathcal{O}_k \cap \mathfrak{g}_k$ is connected.
\end{lemma}
\begin{proof}
Suppose that there exist two nonempty proper open subsets
$U_1,U_2$ of $\mathcal{O}_k \cap \mathfrak{g}_k$ such that
$$\mathcal{O}_k \cap \mathfrak{g}_k = U_1 \cup U_2 \quad
\text{and} \quad U_1\cap U_2 = \varnothing.$$ Then by definition,
we can find two nonempty open subsets $\Omega_1, \Omega_2$ of
$\mathfrak{g}_k$ such that
$$U_i=\Omega_i\cap\mathcal{O}_k \cap \mathfrak{g}_k, \,i=1,2.$$
Since $\Omega_1,\Omega_2$ and $\mathcal{O}_k \cap \mathfrak{g}_k$
are nonempty open subsets of $\mathfrak{g}_k$, all of them are
dense in $\mathfrak{g}_k$. Therefore, $$U_1 \cap U_2 = \Omega_1
\cap \Omega_2 \cap \mathcal{O}_k \cap \mathfrak{g}_k \neq
\varnothing,$$ which contradicts our assumption.
\end{proof}

By Lemmas \ref{lemma:7.4} and \ref{lemma:7.5}, $\mathcal{O}_k \cap
\mathfrak{g}_k$ is an open dense connected subset of
$\mathfrak{g}_k$. Referring to the proof of Theorem \ref{thm:7.1},
we see that $\mathcal{O}_k \cap \mathfrak{g}_k$ is a smooth
manifold and the $G_0$-orbits of $\mathfrak{g}_k$ inside
$\mathcal{O}_k \cap \mathfrak{g}_k$ are precisely its connected
components, which must be $\mathcal{O}_k \cap \mathfrak{g}_k$
itself. Thus $\mathcal{O}_k \cap \mathfrak{g}_k$ is an open dense
$G_0$-orbit in $\mathfrak{g}_k$
\end{proof}

In fact, from Proposition \ref{prop:5.2}, we know that
prehomogeneity is an infinitesimal notion determined only by the
action of the Lie algebra $\mathfrak{g}_0$, so that it depends
only on the Lie type of the complex semisimple Lie group $G$ and
that of the reductive subalgebra $G_0$. This suggests a reason why
this notion should be related to augmentations of Dynkin diagrams,
which capture exactly the Lie types.

Finally, we close this section with a result concerning the corresponding action of the real forms of $G_0$.

\begin{theorem}
Let $(G_0)_\mathbb{R}$ be a real form of $G_0$. Regarding $\mathfrak{g}_i\,(i\neq 0)$ as a complex representation of $(G_0)_\mathbb{R}$ if it is of real type, then the $(G_0)_\mathbb{R}$ action on $(\mathfrak{g}_i)_\mathbb{R}$ has a finite number of orbits. In particular, if $(G_0)_\mathbb{R}$ is the split form of $G_0$, then $\mathfrak{g}_i\,(i \neq 0)$ is always of real type and the corresponding real representation $(\mathfrak{g}_i)_\mathbb{R}$ consists of a finite number of orbits.
\end{theorem}
\begin{proof}
It is a direct consequence of Theorem \ref{thm:7.1} and Theorem \ref{thm:4.6}.
\end{proof}

\section{Termination of $\mathbb{Z}$-Gradings} \label{sect:8}
Recall that upon choosing a simple root $\alpha_0 \in \Pi$, we
have constructed a $\mathbb{Z}$-gradation $$\mathfrak{g}=
\bigoplus_{i \in \mathbb{Z}} \mathfrak{g}_i.$$ Since
$\mathfrak{g}$ is finite dimensional, there exists a positive
integer $n$ such that $\mathfrak{g}_i=0$ for all $|i|>n$. In
virtue of Proposition \ref{prop:3.1}, $\mathfrak{g}_{-i}$ is
naturally identified with $\mathfrak{g}_i^*$ as a
$\mathfrak{g}_0$-representation using the Killing form of
$\mathfrak{g}$. In other words, we have $$\mathfrak{g}=
\bigoplus_{i =-n}^n \mathfrak{g}_i,$$ where $\dim \mathfrak{g}_n =
\dim \mathfrak{g}_{-n} \neq 0$. We call $n$ the order of
$\mathfrak{g}$ with respect to $\alpha_0$ or the order of the
$\mathbb{Z}$-gradation. In fact, there is an easy algorithm to
compute the order $n$. We will consider only the case in which
$\mathfrak{g}$ is simple, the general case follows by considering
the simple ideal containing the root space
$\mathfrak{g}_{\alpha_0}$.

From now on, suppose that $\mathfrak{g}$ is a complex simple Lie algebra. Let $\widetilde{\alpha} \in \Delta$ be the highest root of $\mathfrak{g}$. Then $$\widetilde{\alpha}= \sum_{i=0}^\ell n_i \alpha_i$$ for some positive integers $n_i, i=0,\ldots, \ell$.

\begin{proposition} \label{prop:8.1}
With the above notations, the order of $\mathfrak{g}$ with respect to $\alpha_0$ is $n_0$.
\end{proposition}
\begin{proof}
By definition, we have $\widetilde{\alpha} \in \Delta_{n_0}$ and that for every root $$\alpha= \sum_{i=0}^\ell m_i \alpha_i \in \Delta$$ where $m_i$ are non-negative integers for $i=0,\ldots,\ell$, we have $m_i \leq n_i$ for all $i=0, \ldots, \ell$. In particular, $n_0$ is the greatest integer n for which $\Delta_n \neq \varnothing$.
\end{proof}

Indeed, in each simple case, we can write down the highest root
explicitly. Table \ref{table:4} shows the Dynkin diagrams of all
simple complex Lie algebras with each node labelled by the
coefficient of the corresponding simple root in the highest root,
which is just the order with respect to the corresponding simple
root according to Proposition \ref{prop:8.1}.

\begin{table} \caption{Dynkin diagrams with nodes labelled by the orders with respect to the corresponding simple roots.} \label{table:4}
\begin{displaymath}
\xymatrix{
A_n & \circ \ar@{-}[r]_{\!\!\!\!\!\!\!\!\!\!\!\!\!\!\!\!\!\!\!1} &\circ \ar@{-}[r]_{\!\!\!\!\!\!\!\!\!\!\!\!\!\!\!\!\!\!\!1} &\cdots \ar@{-}[r] &\circ \ar@{-}[r]_{\!\!\!\!\!\!\!\!\!\!\!\!\!\!\!\!\!\!\!1} &\circ \ar@{-}[l]^{\,\,\,\qquad \qquad 1}\\
B_n & \circ \ar@{-}[r]_{\!\!\!\!\!\!\!\!\!\!\!\!\!\!\!\!\!\!\!1} & \circ \ar@{-}[r]_{\!\!\!\!\!\!\!\!\!\!\!\!\!\!\!\!\!\!\!\!2} & \cdots \ar@{-}[r] &\circ \ar@2{->}[r]_{\,\,\,\,\,\,\,2\,\,\,\,\,\,\,\,\qquad 2} &\circ\\
C_n & \circ \ar@{-}[r]_{\!\!\!\!\!\!\!\!\!\!\!\!\!\!\!\!\!\!\!\!2} & \circ \ar@{-}[r]_{\!\!\!\!\!\!\!\!\!\!\!\!\!\!\!\!\!\!\!\!2} & \cdots \ar@{-}[r]_{\,\,\,\,\,\qquad \qquad 2} &\circ &\circ \ar@2{->}[l]^{\,\,\,\qquad \qquad 1}\\
D_n & \circ \ar@{-}[r]_{\!\!\!\!\!\!\!\!\!\!\!\!\!\!\!\!\!\!\!1} &\circ \ar@{-}[r]_{\!\!\!\!\!\!\!\!\!\!\!\!\!\!\!\!\!\!\!\!2} &\cdots \ar@{-}[r] &\circ \ar@{-}[r]_{\!\!\!\!\!\!\!\!\!\!\!\!\!\!\!\!\!\!\!\!2} &\circ \ar@{-}[l]^{\,\,\,\qquad \qquad 1}\\
&&&&\circ \ar@{-}[u] \ar@{-}_{1}\\
E_6 & \circ \ar@{-}[r]_{\!\!\!\!\!\!\!\!\!\!\!\!\!\!\!\!\!\!\!1} &\circ \ar@{-}[r]_{\!\!\!\!\!\!\!\!\!\!\!\!\!\!\!\!\!\!\!\!2} &\circ \ar@{-}[r]_{\!\!\!\!\!\!\!\!\!\!\!\!\!\!\!\!\!\!\!\!3} &\circ \ar@{-}[r]_{\!\!\!\!\!\!\!\!\!\!\!\!\!\!\!\!\!\!\!\!2} &\circ \ar@{-}[l]^{\,\,\,\qquad \qquad 1}\\
&&&\circ \ar@{-}[u] \ar@{-}_{2}\\
E_7 & \circ \ar@{-}[r]_{\!\!\!\!\!\!\!\!\!\!\!\!\!\!\!\!\!\!\!\!2} &\circ \ar@{-}[r]_{\!\!\!\!\!\!\!\!\!\!\!\!\!\!\!\!\!\!\!\!3} &\circ \ar@{-}[r]_{\!\!\!\!\!\!\!\!\!\!\!\!\!\!\!\!\!\!\!\!4} &\circ \ar@{-}[r]_{\!\!\!\!\!\!\!\!\!\!\!\!\!\!\!\!\!\!\!\!3} &\circ \ar@{-}[r]_{\!\!\!\!\!\!\!\!\!\!\!\!\!\!\!\!\!\!\!\!2} &\circ \ar@{-}[l]^{\,\,\,\qquad \qquad 1}\\
&&&\circ \ar@{-}[u] \ar@{-}_{2}\\
E_8 & \circ \ar@{-}[r]_{\!\!\!\!\!\!\!\!\!\!\!\!\!\!\!\!\!\!\!\!2} &\circ \ar@{-}[r]_{\!\!\!\!\!\!\!\!\!\!\!\!\!\!\!\!\!\!\!\!4} &\circ \ar@{-}[r]_{\!\!\!\!\!\!\!\!\!\!\!\!\!\!\!\!\!\!\!\!6} &\circ \ar@{-}[r]_{\!\!\!\!\!\!\!\!\!\!\!\!\!\!\!\!\!\!\!\!5} &\circ \ar@{-}[r]_{\!\!\!\!\!\!\!\!\!\!\!\!\!\!\!\!\!\!\!\!4} &\circ \ar@{-}[r]_{\!\!\!\!\!\!\!\!\!\!\!\!\!\!\!\!\!\!\!\!3} &\circ \ar@{-}[l]^{\,\,\,\qquad \qquad 2}\\
&&&\circ \ar@{-}[u] \ar@{-}_{3}\\
F_4 & \circ \ar@{-}[r]_{\!\!\!\!\!\!\!\!\!\!\!\!\!\!\!\!\!\!\!\!2} &\circ \ar@2{->}[r]_{\!\!\!\!\!\!\!\!\!\!\!\!\!\!\!\!\!\!\!\!\!\!3} &\circ \ar@{-}[r]_{\!\!\!\!\!\!\!\!\!\!\!\!\!\!\!\!\!\!\!\!4} &\circ \ar@{-}[l]^{\,\,\,\,\qquad \qquad 2}\\
G_2 & \circ \ar@3{<-}[r]_{\!\!\!\!\!\!\!\!\!\!\!\!\!\!\!\!\!\!\!\!3} &\circ \ar@{-}_{\,\,\,\,2}
}
\end{displaymath}
\end{table}

From Table \ref{table:4}, we see immediately that the only
possible orders are $1\leq n_0 \leq 6$. Indeed, for those
$\mathfrak{g}$ with order $n_0 >1$ with respect to $\alpha_0$, we
can find a semisimple regular subalgebra
$\widetilde{\mathfrak{g}}$ of $\mathfrak{g}$ containing
$\mathfrak{g}_0$ with the same rank as $\mathfrak{g}$ such that
the corresponding $\mathbb{Z}$-gradation has order $1$ and that
the original $\mathbb{Z}$-gradation factors as a
$\mathbb{Z}_{n_0}$-gradation over $\widetilde{\mathfrak{g}}$.

Note that given a system of simple roots $\Pi=\{\alpha_0, \ldots, \alpha_\ell\}$, we have an extended system of simple roots $\widetilde{\Pi}=\Pi \cup \{-\widetilde{\alpha}\}$ by adjoining the lowest root $-\widetilde{\alpha}$ to it. Set $\Pi_i = \widetilde{\Pi} \backslash \{\alpha_i\}$ and $\Delta(i)= \mathbb{Z}\left\langle \Pi_i\right\rangle \cap \Delta$ for $i = 0,\ldots, \ell$. Then it is known that $\Delta(i)$ forms a reduced root system of $\Delta$ which corresponds to a semisimple subalgebra of $\mathfrak{g}$ of the same rank with a system of simple roots given by $\Pi_i$. Besides, we have the following result concerning the maximal regular reductive subalgebras of $\mathfrak{g}$ which is a direct consequence of a result by Borel-de Siebenthal \cite{ref1} on the maximal closed subroot systems:

\begin{theorem} \label{thm:8.2}
Let $\mathfrak{g}$ be the simple Lie algebra defined above and let $\displaystyle{\widetilde{\alpha}= \!\!\sum_{i=0}^\ell n_i \alpha_i}$ be the highest root with respect to the simple root system $\Pi=\{\alpha_0, \ldots, \alpha_\ell\}$. Then all maximal regular reductive subalgebras can be obtained in one of the following ways:
\begin{enumerate}
    \item when $n_i$ is a prime number, the regular semisimple subalgebra is $$\mathfrak{h} \oplus \bigoplus_{\alpha \in \Delta(i)} \mathfrak{g}_\alpha$$ with root system given by $\Delta(i)$;
    \item when $n_i=1$, the regular reductive subalgebra is $$\mathfrak{h} \oplus \bigoplus_{\alpha \in \Delta_0(i)} \mathfrak{g}_\alpha,$$ where $\Delta_0(i) = \mathbb{Z}\left\langle \Pi \backslash \{\alpha_i\}\right\rangle \cap \Delta$.
\end{enumerate}
\end{theorem}

For a detailed proof see Goto and Grosshans \cite{ref6}. The main idea is that every maximal subroot system is generated by an element $\lambda \in \mathfrak{h}^*$ in the sense of $\{\alpha \in \Delta| \left\langle \lambda, \alpha \right\rangle \in \mathbb{Z}\}$. But the choice of such $\lambda$ is invariant under the affine Weyl group $W_\text{aff}$, which can be assumed to lie in the closure of the fundamental alcove. Finally by explicit case-by-case computations, we obtain the above result.

There is a useful criterion for a regular subalgebra being reductive:

\begin{proposition} \label{prop:8.3}
Let $\mathfrak{f}= \mathfrak{k} \oplus \bigoplus_{\alpha \in \Phi} \mathfrak{g}_\alpha$ be a regular subalgebra of a semisimple Lie algebra $\mathfrak{g}$ with respect to a Cartan subalgebra $\mathfrak{h}$, where $\mathfrak{k}$ is a subspace of $\mathfrak{h}$ and $\Phi\subset \Delta$ as a subroot system. Then $\mathfrak{f}$ is reductive if and only if $\Phi$ is closed and symmetric (i.e. $(\Phi + \Phi)\cap \Delta \subset \Phi$ and $\Phi=-\Phi$) and $\text{span}\{H_\alpha| \alpha \in \Phi\} \subset \mathfrak{k}$.
\end{proposition}

\begin{proposition} \label{prop:8.4}
Let $\displaystyle{\mathfrak{g}= \bigoplus_{i=-n}^{n} \mathfrak{g}_i}$ be the $\mathbb{Z}$-gradation defined above with $\mathfrak{g}_{n} \neq 0$ and $m$ a positive integer dividing $n$. Then $$\mathfrak{g}(m)_{0}:= \bigoplus_{i \equiv 0 (\text{mod } m)} \mathfrak{g}_i$$ is a regular semisimple subalgebra of $\mathfrak{g}$ of the same rank containing $\mathfrak{g}_0$. Moreover, $\mathfrak{g}(n)_0$ is a maximal regular reductive subalgebra of $\mathfrak{g}$ whose root system is isomorphic to $\Delta(0)= \mathbb{Z}\left\langle \Pi_0 \right\rangle \cap \Delta$ and its system of simple roots is $\Pi_0$.
\end{proposition}
\begin{proof}
Note that the set of roots in $\Delta$ occurring in
$\mathfrak{g}(m)_0$ is $$\bigcup_{k \in \mathbb{Z}} \Delta_{km},$$
which is closed and symmetric. Thus, by Proposition
\ref{prop:8.3}, $\mathfrak{g}(m)_0$ is a regular reductive
subalgebra of $\mathfrak{g}$. As $\mathfrak{g}_0 \subset
\mathfrak{g}(m)_0$, $\mathfrak{z}_{\mathfrak{g}(m)_0} \subset
\mathfrak{z}_{\mathfrak{g}_0}=\mathbb{C}c$. But $c$ acts
nontrivially on $\mathfrak{g}_m$ implies that
$\mathfrak{z}_{\mathfrak{g}(m)_0}=0$. Hence $\mathfrak{g}(m)_0$ is
a semisimple subalgebra of the same rank.

To verify the second statement is equivalent to show that $$\mathfrak{g}(n)_0 = \mathfrak{h} \oplus \bigoplus_{\alpha \in \Delta(0)} \mathfrak{g}_\alpha,$$ which reduces to verify that $\Delta(0)= \Delta_{-n} \cup \Delta_0 \cup \Delta_n$. Note that $\widetilde{\alpha} \in \Delta_n$, we have $\Delta_n = (\mathbb{Z}\left\langle \alpha_1, \ldots, \alpha_\ell\right\rangle + \widetilde{\alpha})\cap \Delta$; similarly, we have $\Delta_{-n} = (\mathbb{Z}\left\langle \alpha_1, \ldots, \alpha_\ell\right\rangle - \widetilde{\alpha})\cap \Delta$. Therefore,
\begin{align*}
\Delta(0) &=\mathbb{Z}\left\langle \alpha_1, \ldots, \alpha_\ell, \widetilde{\alpha}\right\rangle \cap \Delta\\
&= \bigcup_{i=-1}^1((\mathbb{Z}\left\langle \alpha_1, \ldots, \alpha_\ell\right\rangle + i\widetilde{\alpha})\cap \Delta)\\
&= \Delta_{-n} \cup \Delta_0 \cup \Delta_n.
\end{align*}
\end{proof}

We have the following characterization of the $\mathbb{Z}$-gradation $\displaystyle{\mathfrak{g}= \bigoplus_{i=-n}^{n} \mathfrak{g}_i}$ when $n=1$ or $n$ is a prime number.

\begin{theorem}
Let $\displaystyle{\mathfrak{g}= \bigoplus_{i=-n}^{n} \mathfrak{g}_i}$ be a $\mathbb{Z}$-gradation of $\mathfrak{g}$ as constructed above with $\mathfrak{g}_n \neq 0$.
\begin{enumerate}
    \item If $n=1$, then $\mathfrak{g}_0$ is a maximal reductive subalgebra of $\mathfrak{g}$ with an one-dimensional center;
    \item if $n$ is a prime number, then $\mathfrak{g}(n)_0$ is a maximal semisimple subalgebra of $\mathfrak{g}$ of the same rank for which $\mathfrak{g}_0$ lies in $\mathfrak{g}(n)_0$ as a maximal reductive subalgebra with a one-dimensional center.
\end{enumerate}
\end{theorem}
\begin{proof}
By Proposition \ref{prop:8.4}, we have $\mathfrak{g}_0$ in (i) being the second case of Theorem \ref{thm:8.2} and $\mathfrak{g}(n)_0$ in (ii) being the first case of Theorem \ref{thm:8.2}.
\end{proof}

\begin{proposition}
Under the same conditions as in Proposition \ref{prop:8.4}, we have $$\mathfrak{g} = \bigoplus_{j \in \mathbb{Z}_m} \mathfrak{g}(m)_j$$ is a $\mathbb{Z}_m$-gradation of $\mathfrak{g}$, where $\displaystyle{\mathfrak{g}(m)_j= \bigoplus_{i \equiv j (\text{mod }m)}\mathfrak{g}_i}$.
\end{proposition}
\begin{proof}
For $j_1, j_2 \in \mathbb{Z}_m$ and $i_s \equiv j_s (\text{mod } m)$, $s=1,2$, we have $i_1+i_2 \equiv j_i + j_2 (\text{mod } m)$ and $[\mathfrak{g}_{i_1}, \mathfrak{g}_{i_2}] \subset \mathfrak{g}_{i_1 + i_2} \subset \mathfrak{g}(m)_{j_1 + j_2}$. Hence $$[\mathfrak{g}(m)_{j_1}, \mathfrak{g}(m)_{j_2}]\subset \mathfrak{g}(m)_{j_1+j_2}.$$
\end{proof}

\section{Explicit Construction of Generic Elements in Simply-laced Cases} \label{sect:9}
Throughout this section, $\mathfrak{g}$ is assumed to be simply-laced, i.e. the corresponding Dynkin diagram $\Gamma(\mathfrak{g})$ consists of single edges only. Assume that the nodes of $\Gamma(\mathfrak{g})$ are indexed by a system of simple roots $\Pi=\{\alpha_0, \ldots, \alpha_\ell\}$. In this case, the corresponding Cartan matrix $(\left\langle \alpha_i,\alpha_j\right\rangle_\mathfrak{h})_{i,j =0,\ldots, \ell}$ is completely determined by their restrictions onto $\mathfrak{t}$, namely $\left\langle \alpha_i,\alpha_j\right\rangle_\mathfrak{h} = \left\langle \alpha_i|_\mathfrak{t},\alpha_j|_\mathfrak{t}\right\rangle_\mathfrak{t}$ for $j\neq 0, i=0, \ldots, \ell$, and $\left\langle \alpha_k,\alpha_0\right\rangle_\mathfrak{h}=\left\langle \alpha_0, \alpha_k\right\rangle_\mathfrak{h} = \left\langle \alpha_0|_\mathfrak{t},\alpha_k|_\mathfrak{t}\right\rangle_\mathfrak{t}$ according to Lemma \ref{lemma:6.1} and the fact that all roots have the same length in a simply-laced semisimple Lie algebra.

Now we denote $W_0$ as the Weyl group of $\mathfrak{g}_0^{ss}$
generated by the reflections in $\mathfrak{t}^*$ along
$\{\alpha_1|_\mathfrak{t}, \ldots, \alpha_\ell|_\mathfrak{t}\}$.
Let $\widetilde{W_0}$ be the subgroup of the Weyl group $W$ of
$\mathfrak{g}$ generated by reflections along $\{\alpha_1, \ldots,
\alpha_\ell\}$. As every element in $\widetilde{W_0}$ preserves
the subspace $\mathfrak{t}^*$, the natural restriction map induces
an isomorphism from $\widetilde{W_0}$ to $W_0$, denoted by
$\widetilde{w}\mapsto w$. Also, we see that $\widetilde{W_0}$
stabilizes each $\Delta_i$.

Let $(H_\alpha,X_\alpha, Y_\alpha)$ be a standard $\mathfrak{sl}_2$-triple corresponding to $\alpha \in \Delta^+$, which will be fixed once and for all throughout the whole section. In the following, we will give an explicit construction of generic elements in $\mathfrak{g}_1$ and $\mathfrak{g}_{-1}$ as $G_0$ representations.

\begin{lemma} \label{lemma:9.1}
$\widetilde{w}(\alpha)|_\mathfrak{t}=w(\alpha|_\mathfrak{t})$ for all $\widetilde{w} \in \widetilde{W_0}, \alpha \in \Delta$. In particular, we have $|\widetilde{W_0}\cdot \alpha_0| = \left|W_0 \cdot \alpha_0|_\mathfrak{t}\right|$.
\end{lemma}
\begin{proof}
For all $j=1,\ldots, \ell$, by Lemma \ref{lemma:6.1}, we have
\begin{align*}
(\widetilde{w}(\alpha)|_\mathfrak{t},\alpha_j|_\mathfrak{t})_\mathfrak{t} &= (\widetilde{w}(\alpha), \alpha_j)_\mathfrak{h}\\
&=(\alpha, \widetilde{w}^{-1}(\alpha_j))_\mathfrak{h}\\
&=(\alpha|_\mathfrak{t},\widetilde{w}^{-1}(\alpha_j)|_\mathfrak{t})_\mathfrak{t}\\
&=(\alpha|_\mathfrak{t},w^{-1}(\alpha_j|_\mathfrak{t}))_\mathfrak{t}\\
&=(w(\alpha|_\mathfrak{t}), \alpha_j|_\mathfrak{t})_\mathfrak{t}.
\end{align*}
Since $\{\alpha_1, \ldots, \alpha_\ell \}$ form a basis of $\mathfrak{t}^*$, we conclude that $\widetilde{w}(\alpha)|_\mathfrak{t}=w(\alpha|_\mathfrak{t})$.
\end{proof}

\begin{lemma} \label{lemma:9.2}
For $k \neq 0$, $$\sum_{\alpha \in \widetilde{W_0} \cdot \alpha_0} \left\langle \alpha,\alpha_k \right\rangle_\mathfrak{h} =0.$$
\end{lemma}
\begin{proof}
Note that
\begin{align*}
\sum_{\alpha \in \widetilde{W_0}\cdot \alpha_0} \left\langle \alpha,\alpha_k\right\rangle_\mathfrak{h} &= \sum_{\alpha \in \widetilde{W_0}\cdot \alpha_0} \left\langle \alpha|_\mathfrak{t},\alpha_k|_\mathfrak{t} \right\rangle_\mathfrak{t} &(\text{Lemma \ref{lemma:6.1}})\\
&= \sum_{\beta \in W_0 \cdot \alpha_0|_\mathfrak{t}} \left\langle \beta,\alpha_k|_\mathfrak{t} \right\rangle_\mathfrak{t} &(\text{Lemma \ref{lemma:9.1}})\\
&= \left\langle \sum_{\beta \in W_0 \cdot \alpha_0|_\mathfrak{t}} \beta,\alpha_k|_\mathfrak{t} \right\rangle_\mathfrak{t}.
\end{align*}
Here $\displaystyle{\sum_{\beta \in W_0 \cdot \alpha_0|_\mathfrak{t}} \beta}$ is just the sum of all extremal weights of $\mathfrak{g}_1$ as a $\mathfrak{g}_0^{ss}$ representation. As the set of extremal weights is symmetric about the origin, we have $\displaystyle{\sum_{\beta \in W_0 \cdot \alpha_0|_\mathfrak{t}} \beta=0}$ and thus $$\sum_{\alpha \in \widetilde{W_0} \cdot \alpha_0} \left\langle \alpha,\alpha_k \right\rangle_\mathfrak{h} =0.$$
\end{proof}

\begin{lemma} \label{lemma:9.3}
$$\sum_{\alpha \in \widetilde{W_0}\cdot \alpha_0} \left\langle \alpha,\alpha_0\right\rangle_\mathfrak{h}= 2 \left|W_0 \cdot \alpha_0|_\mathfrak{t} \right| \cdot \frac{\left\|\alpha_0'\right\|^2}{\left\|\alpha_0\right\|^2},$$ where $\alpha_0'$ is the orthogonal projection of $\alpha_0$ to $(\mathbb{C}c)^*$, i.e. $\alpha_0 = \alpha_0|_\mathfrak{t} + \alpha_0'$ with $\left\langle \alpha|_\mathfrak{t},\alpha_0' \right\rangle_\mathfrak{h}=0$.
\end{lemma}
\begin{proof}
First note that $\widetilde{W_0}\cdot \alpha_0 \subset \Delta_1$ so that every element in $\widetilde{W_0}\cdot \alpha_0$ has the same orthogonal projection $\alpha_0'$ onto $(\mathbb{C}c)^*$. Then by applying the above two lemmas, we obtain
\begin{align*}
\sum_{\alpha \in \widetilde{W_0} \cdot \alpha_0} \alpha &= \sum_{\alpha \in \widetilde{W_0} \cdot \alpha_0} (\alpha|_\mathfrak{t} + \alpha_0')\\
&= \sum_{\beta \in W_0 \cdot \alpha_0|_\mathfrak{t}} \beta + \sum_{\alpha \in \widetilde{W_0} \cdot \alpha_0} \alpha_0'\\
&= \left|\widetilde{W_0} \cdot \alpha_0 \right|\cdot \alpha_0'\\
&= \left|W_0 \cdot \alpha_0|_\mathfrak{t} \right|\cdot \alpha_0'.
\end{align*}
It follows that
\begin{align*}
\sum_{\alpha \in \widetilde{W_0}\cdot \alpha_0} \left\langle \alpha,\alpha_0\right\rangle_\mathfrak{h} &= \left|W_0 \cdot \alpha_0|_\mathfrak{t} \right| \left\langle \alpha_0',\alpha_0\right\rangle_\mathfrak{h}\\
&= 2 \left|W_0 \cdot \alpha_0|_\mathfrak{t} \right| \cdot \frac{\left( \alpha_0', \alpha_0|_\mathfrak{t} + \alpha_0'\right)}{\left(\alpha_0,\alpha_0\right)}\\
&= 2 \left|W_0 \cdot \alpha_0|_\mathfrak{t} \right| \cdot \frac{\left\|\alpha_0'\right\|^2}{\left\|\alpha_0\right\|^2}.
\end{align*}
\end{proof}

\begin{theorem}
Let $\mathfrak{g}$ be a simply-laced semisimple Lie algebra as
defined above. If $(G_0,\mathfrak{g}_{-1})$ is a regular prehomogeneous space, then there exist $X \in \mathfrak{g}_1, Y \in
\mathfrak{g}_{-1}$ such that $[X,Y]=c$.
\end{theorem}
\begin{proof}
For $k = 0,\ldots, \ell$, we have
\begin{align*}
&\left(H_{\alpha_k}, \left[ \sum_{\alpha \in \widetilde{W_0} \cdot \alpha_0} X_\alpha, \sum_{\beta \in \widetilde{W_0} \cdot \alpha_0} Y_\beta \right]\right)\\
=& \left(\left[H_{\alpha_k}, \sum_{\alpha \in \widetilde{W_0} \cdot \alpha_0} X_\alpha \right], \sum_{\beta \in \widetilde{W_0} \cdot \alpha_0} Y_\beta \right)\\
=&\left(\sum_{\alpha \in \widetilde{W_0} \cdot \alpha_0} \alpha(H_{\alpha_k}) X_\alpha, \sum_{\alpha \in \widetilde{W_0} \cdot \alpha_0} Y_\alpha \right)\\
=&\sum_{\alpha \in \widetilde{W_0} \cdot \alpha_0} \left\langle \alpha,\alpha_k\right\rangle_\mathfrak{h} \left(X_\alpha, \sum_{\alpha \in \widetilde{W_0} \cdot \alpha_0} Y_\alpha\right)\\
=&\sum_{\alpha \in \widetilde{W_0} \cdot \alpha_0} \left\langle \alpha,\alpha_k\right\rangle_\mathfrak{h} (X_\alpha, Y_\alpha)\\
=& \frac{1}{2} \sum_{\alpha \in \widetilde{W_0} \cdot \alpha_0} \left\langle \alpha, \alpha_k\right\rangle_\mathfrak{h} (H_\alpha, H_\alpha)\\
=& \frac{1}{2} \sum_{\alpha \in \widetilde{W_0} \cdot \alpha_0} \left\langle \alpha,\alpha_k\right\rangle_\mathfrak{h} (H_{\alpha_0}, H_{\alpha_0})\\
=&
\begin{cases}
0, &\text{ if } k \neq 0\,(\text{Lemma \ref{lemma:9.2}}),\\
\left|W_0 \cdot \alpha_0|_\mathfrak{t}\right| \cdot \frac{\left\|\alpha_0'\right\|^2}{\left\|\alpha_0\right\|^2} \cdot (H_{\alpha_0}, H_{\alpha_0}), &\text{ if } k=0\,(\text{Lemma \ref{lemma:9.3}}).
\end{cases}
\end{align*}
Set $$X= \frac{1}{\sqrt{2 \left|W_0 \cdot \alpha_0|_\mathfrak{t}\right|}} \cdot \frac{\left\|\alpha_0\right\|}{\left\|\alpha_0'\right\|} \sum_{\alpha \in \widetilde{W_0} \cdot \alpha_0} X_\alpha, Y= \frac{1}{\sqrt{2 \left|W_0 \cdot \alpha_0|_\mathfrak{t}\right|}} \cdot \frac{\left\|\alpha_0\right\|}{\left\|\alpha_0'\right\|} \sum_{\alpha \in \widetilde{W_0} \cdot \alpha_0} Y_\alpha.$$ Then $X \in \mathfrak{g}_1$, $Y \in \mathfrak{g}_{-1}$ and that
\begin{align*}
\alpha_k([X,Y]) &= \frac{2 \left( H_{\alpha_k},[X,Y]\right)}{(H_{\alpha_k}, H_{\alpha_k})}\\
&=
\begin{cases}
0, \text{ if }k \neq 0,\\
1, \text{ if }k=0.
\end{cases}
\end{align*}
Note that the regularity condition is to ensure that $[X,Y]$ lies in $\mathfrak{h}$ as in this case $$\left[ \sum_{\alpha \in \widetilde{W_0} \cdot \alpha_0} X_\alpha, \sum_{\beta \in \widetilde{W_0} \cdot \alpha_0} Y_\beta \right] = \sum_{\alpha \in \widetilde{W_0} \cdot \alpha_0} \left[ X_\alpha, Y_\alpha \right].$$ It follows that $[X,Y]=c$.
\end{proof}

\begin{corollary}\footnote{This corollary is first proved in \cite{ref11}.}
Let $\mathfrak{g}$ be a simply-laced semisimple Lie algebra as
defined above and $X \in \mathfrak{g}_1$. If there exists $Y \in
\mathfrak{g}_{-1}$ such that $[X,Y]=c$, then
$[X,\mathfrak{g}_0]=\mathfrak{g}_1$, i.e. $X$ is a generic element
of $(G_0,\mathfrak{g}_1)$
\end{corollary}
\begin{proof}
First note that $\frac{1}{2}c,X,Y$ form a standard set of
generators for an $\mathfrak{sl}_2$ subalgebra $\mathfrak{a}$ of
$\mathfrak{g}$. Then from $\mathfrak{sl}_2$ theory, if we
decompose $\mathfrak{g}$ into irreducible $\mathfrak{a}$
representations, there are no weight spaces with weight $1$ and
that $\mathfrak{g}_1$ is the direct sum of all weight spaces of
weight $2$. Hence, we have $[X,\mathfrak{g}_0]=\mathfrak{g}_1$.
The stabilizer $(\mathfrak{g}_0)_X$ is reductive as it is the
centralizer of $\mathfrak{a}$.
\end{proof}
\section{The Ambient Lie Algebras of Parabolic PVS's} \label{sect:10}
With the effort of the previous sections, we can already conclude
Theorem \ref{thm:1.2}(iii) on the class of simply-laced Lie
algebras. The general situation is more complicated as there can
be more than one ambient Lie algebras $\mathfrak{g}$ associated to
an irreducible prehomogeneous vector space $(G_0,V)$. For
instance, if we consider the prehomogeneous vector space $(GL_2,
\mathbb{C}^2)$, we can choose $\mathfrak{g}$ to be either
$\mathfrak{sl}_3, \mathfrak{so}_5, G_2$.

In this section, we will finish the proof of Theorem \ref{thm:1.2}
by showing that given an irreducible parabolic PVS $(G_0, V, \nu)$
with connecting multiplicities there exists exactly one ambient
Lie algebra $\mathfrak{g}$ containing $\mathfrak{g}_0$ for which
$(\Gamma(\mathfrak{g}),\Gamma(\mathfrak{g}_0^{ss}))$ maps to
$(G_0, V, \nu(\mathfrak{g},V))$ under the correspondence set up in
Theorem \ref{thm:1.2}(iv). The main result we used here is the
Serre's Theorem which states that given a Cartan matrix
$A=(a_{ij})$ of rank $\ell$ there is a semisimple Lie algebra with
$3 \ell$ generators $\{H_i,X_i,Y_i\}, i=1, \ldots , \ell$,
satisfying
\begin{align}\label{eq:10.1}
&[H_i, H_j]=0\\
&[X_i,Y_i]=H_i, [X_i,Y_j] =0 \text{ if } i \neq j\\
&[H_i, X_j]= a_{ji} X_j, [H_i,Y_j]= -a_{ji} Y_j\\
&(ad X_i)^{- a_{ji}+1}(X_j)=0\\  \label{eq:10.5}
&(ad Y_i)^{-a_{ji}+1}(Y_j)=0
\end{align}
unique up to isomorphism.

Suppose $\mathfrak{g}_0^{ss}$ is of rank $\ell$ and by choosing a Cartan subalgebra $\mathfrak{t}$ as usual, we obtain a corresponding root system $\Delta_0$. Finally, we fix a system of simple roots $\Pi_0=\{\alpha_1, \ldots, \alpha_\ell\}$ of $\Delta_0$. To each simple root $\alpha_i$, we already have $\{H_i,X_i,Y_i\}$ satisfying relations in (\ref{eq:10.1})-(\ref{eq:10.5}). Let $\mathfrak{h}=\mathfrak{z}_{\mathfrak{g}_0} \oplus \mathfrak{t}$. To construct $\mathfrak{g}$ it suffices to find $H_0,X_0,Y_0$ which are compatible with other $H_i,X_i,Y_i$.

Let $\nu=(a_1, \ldots, a_r)$ and $\omega$ be the highest weight of
the irreducible representation $(\pi,V)$ of $\mathfrak{g}_0^{ss}$.
Without loss of generality, we can assume that all $a_i \neq 0$
and that the matrix
\begin{align*}
A= \left(%
\begin{array}{ccccccc}
  2 & -\langle\omega, \alpha_1\rangle & \cdots & -\langle\omega, \alpha_r\rangle & 0 & \cdots & 0 \\
  -a_1 & \langle\alpha_1, \alpha_1\rangle & \cdots & \langle\alpha_1, \alpha_r\rangle & \langle\alpha_1, \alpha_{r+1}\rangle & \cdots & \langle\alpha_1, \alpha_\ell\rangle \\
  \vdots &  &  &  &  &  &  \\
  -a_r & \vdots &  &  &  &  & \vdots \\
  0 & \vdots &  &  &  &  & \vdots \\
  \vdots &  &  &  &  &  &  \\
  0 & \langle\alpha_\ell, \alpha_1\rangle &  & \cdots & \cdots &  & \langle\alpha_\ell, \alpha_\ell\rangle \\
\end{array}%
\right)
\end{align*}
is a Cartan matrix, so that we have
$$\frac{a_i(H_i,H_i)}{\omega(H_i)}=K$$ for some fixed nonzero
constant $K$. Then there exists a unique element $H \in
\mathfrak{t}$ such that
\begin{align*}
\alpha_i(H)= \begin{cases}
    -a_i & \text{, if } i=1, \ldots r,\\
    0 & \text{, otherwise.}
\end{cases}
\end{align*}
Pick any nonzero $X_0 \in V^*_{-\omega}$, we can find a unique $c
\in \mathfrak{z}_{\mathfrak{g}_0}$ such that $$\pi^*(H+c)X_0 =
2X_0.$$ Let $\kappa$ be the unique $G_0$-invariant nondegenerate
bilinear form on $\mathfrak{g}_0$ extending the Killing form of
$\mathfrak{g}_0^{ss}$ and satisfies $$\kappa(c,c)= K-(H,H).$$
Note that $\kappa(c, \mathfrak{g}_0^{ss})=0$ is automatic from the
invariance property of $\kappa$, it follows that
$$\kappa(H_0,H_0)= \kappa(H,H)+ \kappa(c,c)=K.$$ Choose $Y_0 \in
V_\omega$ such that $X_0(Y_0)=-\frac{K}{2}$.

Formally, we can now impose conditions
(\ref{eq:10.1})-(\ref{eq:10.5}) to $\{H_i,X_i,Y_i\}_{i=0}^\ell$,
whence we obtain a semisimple Lie algebra $\mathfrak{g}$ by
applying the Serre's theorem. Since the last $3\ell$ generators
$\{H_i,X_i,Y_i\}_{i=1}^\ell$ are also generators for
$\mathfrak{g}_0^{ss}$, we obtain an embedding $\mathfrak{g}_0^{ss}
\subset \mathfrak{g}$, and that $\mathfrak{g}_0 = \mathbb{C}H_0
\oplus \mathfrak{g}_0^{ss} \subset \mathfrak{g}$. The remaining
task is to construct the bracket relations between elements in
$\mathfrak{g}_0,V$ and $V^*$ which coincide with that abstractly
defined in terms of the generators of $\mathfrak{g}$. The obvious
choice of defining the bracket on $\mathfrak{g}_0 \times V$ and
$\mathfrak{g}_0 \times V^*$ is $$[Z,v]= \pi(Z)v,\quad
[Z,f]=\pi^*(Z)f$$ for all $Z \in \mathfrak{g}_0, v \in V, f \in
V^*$. Direct checking shows that (\ref{eq:10.1})-(\ref{eq:10.5})
are satisfied except the equality $[X_0,Y_0]=H_0$ has not yet been
established.

Since $X_0$ is a lowest weight vector of $(\pi^*,V^*)$, for all $i=1, \ldots , \ell$, $$(\pi^*(Y_i)X_0)Y_0 = 0 = \kappa(H_0,Y_i).$$ Similarly, as $Y_0$ is a highest weight vector of $(\pi,V)$, for all $i=1, \ldots \ell$, $$X_0(\pi(X_i))Y_0 = 0 = \kappa(H_0, X_i).$$ Clearly, for $i=r+1, \ldots, \ell$, $a_i=\omega(H_i)=0$, and $$X_0(\pi(X_i))Y_0 = \omega(H_i) X_0(Y_0) = 0 = \kappa(H_0,H_i).$$ For $i=1, \ldots, r$,
\begin{align*}
X_0(\pi(X_i))Y_0 &= \omega(H_i) X_0(Y_0)\\
&= - \omega(H_i)\frac{K}{2}\\
&= - \frac{a_i \kappa(H_i,H_i)}{2}\\
&= \frac{\alpha_i(H)(H_i,H_i)}{2}\\
&= (H,H_i)\\
&= \kappa(H_0,H_i).
\end{align*}
All together, we get
\begin{align} \label{eq:10.6}
\kappa(H_0,Z)=X_0(\pi(Z)Y_0)=-(\pi^*(Z)X_0)Y_0 \qquad \text{for
all } Z \in \mathfrak{g}_0.
\end{align}
From $\pi: \mathfrak{g}_0 \rightarrow End(V)= V^* \otimes V$, we
get the moment map $$\mu_\pi:V^* \times V \rightarrow
\mathfrak{g}_0^*.$$ By identifying $\mathfrak{g}_0$ and
$\mathfrak{g}_0^*$ through $\kappa$, we get a bilinear map
$$\phi:V^* \times V \rightarrow \mathfrak{g}_0.$$ Explicitly,
given $v \in V, f \in V^*$, $\phi(f,v)$ is the unique element such
that $$\kappa(\phi(f,v),Z) = f(\pi(Z)v) = -(\pi^*(Z)f)v \qquad
\text{for all }Z \in \mathfrak{g}_0.$$ In view of (\ref{eq:10.6}),
we have $\phi(X_0,Y_0) = H_0$. Thus this map coincides with the
bracket structure constructed on $\mathfrak{g}$. In other words,
we have $V$ and $V^*$ embedded into $\mathfrak{g}$ with the
bracket between elements of $V$ and $V^*$ given by the map $\phi$.
In particular , we have $X_0$ being a root vector corresponding to
a root $\widetilde{\alpha_0}$ of $\mathfrak{g}$ with respect to
$\mathfrak{h}$. Let $\widetilde{\alpha_i}\in \mathfrak{h}^*$ is
the extension of $\alpha_i \in \mathfrak{t}^*$ by setting
$\widetilde{\alpha_i}(c)=0$. Then we can easily see that $\Pi=
\{\widetilde{\alpha_0}, \widetilde{\alpha_1}, \ldots,
\widetilde{\alpha_\ell}\}$ is a system of simple roots to
$\mathfrak{g}$ with $\widetilde{\alpha_0}|_\mathfrak{t}= -\omega,
\widetilde{\alpha_i}|_\mathfrak{t} = \alpha_i$ for $i=1,\ldots
\ell$, and the corresponding Cartan matrix is given by $A$. This
complete the proof of Theorem \ref{thm:1.2}.
\section{PVS's of Twisted Affine Type} \label{sect:11}
As we have mentioned in the introduction, there are some examples
of prehomogeneous vector spaces consisting of finitely many orbits
which are not of parabolic type. Among the irreducible reduced
ones, there are six exceptional cases as listed in Table
\ref{table:3}. We will briefly explain these structures and will
find a unified way of constructing them.

Observe that we have the following grading:
\begin{align*}
\mathfrak{so}_{10} &= \mathfrak{gl}_1 \times G_2 \oplus \mathbb{C}^7 \oplus (\mathbb{C}^7)^* \oplus \mathbb{C}^7 \oplus (\mathbb{C}^7)^* \oplus \mathbb{C} \oplus (\mathbb{C})^*\\
E_6 &= \mathfrak{gl}_2 \!\times\! G_2 \!\oplus\! (\mathbb{C}^2 \!\otimes \mathbb{C}^7) \!\oplus\! (\mathbb{C}^2 \!\otimes \mathbb{C}^7)^* \!\oplus\! (\mathbb{C}^2 \!\otimes \mathbb{C}^7) \!\oplus\! (\mathbb{C}^2 \!\otimes \mathbb{C}^7)^* \!\oplus\! \mathbb{C}^2 \!\oplus\! (\mathbb{C}^2)^*\\
\mathfrak{so}_{12} &= \mathfrak{gl}_2 \times \mathfrak{so}_7 \oplus \mathbb{C}^7\oplus (\mathbb{C}^2 \otimes S) \oplus (\mathbb{C}^2 \otimes S)^* \oplus \mathbb{C} \oplus \mathbb{C}^*\\
E_7 &= \mathfrak{gl}_3 \times \mathfrak{so}_7 \oplus \mathbb{C}^7 \oplus (\mathbb{C}^3 \otimes S) \oplus (\mathbb{C}^3 \otimes S)^* \oplus (\mathbb{C}^3 \otimes S) \oplus (\mathbb{C}^3 \otimes S)^*\\
E_6 &= \mathfrak{gl}_1 \times \mathfrak{so}_9 \oplus \mathbb{C}^9 \oplus (\mathbb{C} \otimes S) \oplus (\mathbb{C} \otimes S)^*\\
E_7 &= \mathfrak{gl}_1 \times \mathfrak{so}_{11} \oplus \mathbb{C}^{11} \oplus (\mathbb{C} \otimes S) \oplus (\mathbb{C} \otimes S)^* \oplus \mathbb{C} \oplus \mathbb{C}^*
\end{align*}
They are obtained from successive $\mathbb{Z}$-gradations and then
further decomposed by an outer automorphism of the $0^{th}$ graded
reductive subalgebra. For example, $(GL_1 \times G_2,
\mathbb{C}\otimes \mathbb{C}^7)$ can be obtained by first
considering the $\mathbb{Z}$-gradation
$$\mathfrak{so}_{10} = \mathfrak{gl}_1 \times \mathfrak{so}_8
\oplus \mathbb{C}^8 \oplus \left(\mathbb{C}^8\right)^*$$
associated to the augmentation of Dynkin diagrams $(D_5,D_4)$ and
then further decompose the gradation into irreducible
representations of the fixed point subalgebra $G_2$ of
$\mathfrak{so}_8$ by an outer automorphism induced from the
triality of the Dynkin diagram $D_4$. Other cases can be done
similarly by a suitable reduction of their Dynkin diagrams to the
one possessing a nontrivial outer automorphism and then decompose
the gradation by the fixed point subalgebra obtained from the
corresponding outer automorphism. The advantage of doing this is
that the irreducible subrepresentations contained in any nonzero
components of the $\mathbb{Z}$-gradations still lie in the
nilpotent cone of the orginal ambient semisimple Lie algebra, so
that our previous arguments in Section \ref{sect:7} are still
valid in these cases according to Remark 4 after the proof of
Theorem \ref{thm:7.1}. Collectively speaking, the six cases above
can be obtained from the twisted affine diagrams
\begin{displaymath}
\xymatrix{
E_6^{(3)} & \bullet \ar@{-}[r] &\circ \ar@3{<-}[r] &\circ\\
E_6^{(2)} & \circ \ar@{-}[r] &\bullet \ar@{-}[r] &\circ \ar@2{<-}[r] &\circ \ar@{-}[r] &\circ\\
E_7^{(2)}  &\bullet \ar@{-}[r] &\circ \ar@2{<-}[r] &\circ \ar@{-}[r] &\circ \ar@{-}[r] &\circ\\
}
\end{displaymath}
by deleting the painted node. For example, in the case of $(GL_1
\times G_2, \mathbb{C}\otimes \mathbb{C}^7)$, the naive way to
associate the twisted affine diagram is the construct an
augmentation of the Dynkin diagram $G_2$ by adjoining the lowest
weight of $\mathbb{C}\otimes \mathbb{C}^7$ to the corresponding
system of simple roots. It also works for the cases $(GL_2 \times
Spin_7, \mathbb{C}^2 \otimes S)$ and $(GL_1 \times Spin_9,
\mathbb{C}\otimes S)$. In fact, these pieces of information give
rise to $\mathbb{Z}_m$-gradations instead of
$\mathbb{Z}$-gradations since they can be treated as the fixed
point algebra of appropriate outer automorphisms of a regular
subalgebra of the ambient Lie algebra, and then the corresponding
branching of the adjoint representations yields the above
decompositions. These $\mathbb{Z}_m$-gradations possess an extra
symmetry between the graded pieces which allow us to extend the
symmetry group to $GL_2 \times G_2, GL_3 \times Spin_7$ and $GL_1
\times Spin_{11}$ respectively in the remaining three cases.
From this point of view, it is reasonable to call them the
prehomogeneous vector spaces of twisted affine type.

\section{Orbit Structure of $(GL_2 \times SL_{2m+1}, \mathbb{C}^2 \otimes \Lambda^2 \mathbb{C}^{2m+1})$} \label{sect:12}
In this section, we will examine the orbit structure of the
exceptional series $(GL_2 \times SL_{2m+1}, \mathbb{C}^2 \otimes
\Lambda^2 \mathbb{C}^{2m+1}), m \geq 4$, of irreducible reduced
PVS's consisting of an infinite number of orbits. Basically, the
reason of having an infinite number of orbits is due to the
absence of an open orbit in $(GL_2 \times SL_{2m}, \mathbb{C}^2
\otimes \Lambda^2 \mathbb{C}^{2m})$. At the same time, the
construction given below also explains why it is not the case when
$m \leq 3$.

First, we decompose $\mathbb{C}^2 \otimes \Lambda^2
\mathbb{C}^{2m+1}$ into two two parts:
$$\mathbb{C}^2 \otimes
\Lambda^2 \mathbb{C}^{2m+1} = \{(1,0) \otimes \omega_1 + (0,1)
\otimes \omega_2|\omega_1,\omega_2 \in \Lambda^2 \mathbb{C}^{2m+1}
\} = U_1 \cup U_2$$ where
\begin{align*}
U_1 &=\{(1,0) \otimes \omega_1 + (0,1) \otimes \omega_2| \omega_1,
\omega_2 \in \Lambda^2 V \text{ for some $2m$-dim'$\ell$ } V
\subset \mathbb{C}^{2m+1} \},\\
U_2 &= \mathbb{C}^2 \otimes \Lambda^2 \mathbb{C}^{2m+1} - U_1.
\end{align*}
Note that for any $(A,g) \in GL_2 \times SL_{2m+1}$, $A=\left(%
\begin{array}{cc}
  a & b \\
  c & d \\
\end{array}%
\right)$, $(1,0) \otimes \omega_1 + (0,1) \otimes \omega_2 \in
\mathbb{C}^2 \otimes \Lambda^2 \mathbb{C}^{2m+1}$, we have
\begin{align*}
&(A,g)\cdot [(1,0) \otimes \omega_1 + (0,1) \otimes \omega_2]\\
=&(a,c) \otimes g \omega_1 + (b,d) \otimes g \omega_2\\
=&(1,0) \otimes (a g\omega_1 + b g\omega_2) + (0,1) \otimes (c
g\omega_1+ d g\omega_2).
\end{align*}
So if both $\omega_1, \omega_2 \in \Lambda^2 V$ for some $V
\subset \mathbb{C}^{2m+1}$, $g\omega_1, g\omega_2 \in gV$ and
$(1,0) \otimes (a g\omega_1 + b g\omega_2) + (0,1) \otimes (c
g\omega_1+ d g\omega_2) \in U_1$. It follows that $U_1$ and $U_2$
are $GL_2 \times SL_{2m+1}$-invariant subsets.

We fix the standard $\mathbb{C}^{2m}$ as generated by the first
$2m$ coordinate vectors $e^1, \ldots, e^{2m}$ of
$\mathbb{C}^{2m+1}$. We see that any $GL_2 \times SL_{2m+1}$-orbit
in $U_1$ intersects $\Lambda^2\mathbb{C}^{2m}$ nontrivially as a
$GL_2 \times SL_{2m}$-orbit in $\Lambda^2\mathbb{C}^{2m}$. In
other words, we have an one-to-one correspondence between the $GL_2
\times SL_{2m+1}$-orbits in $U_1$ and the $GL_2 \times
SL_{2m}$-orbits in $\Lambda^2\mathbb{C}^{2m}$.

To see that $\Lambda^2\mathbb{C}^{2m}$ has infinitely many $GL_2
\times SL_{2m}$-orbits, we attach to each $(1,0) \otimes \omega_1
+ (0,1) \otimes \omega_2$ a two parameter family of top exterior
forms $$(\lambda \omega_1 + \mu \omega_2)^m =
f_{(\omega_1,\omega_2)}(\lambda,\mu) \, e^1 \wedge \cdots \wedge
e^{2m},$$ where $f_{(\omega_1,\omega_2)}$ is a homogeneous
polynomial of degree $m$ in $\lambda,\mu$. It is easy to check
that such polynomials satisfies
$$f_{A(g\omega_1,g\omega_2)}(\lambda,\mu) =
f_{(\omega_1,\omega_2)}((\lambda,\mu)A)$$ for any $(A,g) \in GL_2
\times SL_{2m}$. Thus we obtain a map
\begin{equation} \label{eqt:12.1}
    \Phi:\mathbb{C}^2 \otimes
\Lambda^2\mathbb{C}^{2m}/ GL_2 \times SL_{2m} \longrightarrow
S^m\mathbb{C}^2/GL_2,
\end{equation}
where $S^m\mathbb{C}^2$ is identified with the space of
homogeneous polynomials of degree $m$.

\begin{proposition} \label{prop:12.1}
The map $\Phi$ defined in $(\ref{eqt:12.1})$ is surjective.
\end{proposition}
\begin{proof}
Pick any nonzero homogeneous polynomial $f(\lambda,\mu)$ of degree
$m$, there exists $p_1=[\lambda_1:\mu_1], \ldots,
p_m=[\lambda_m:\mu_m] \in \mathbb{CP}^1$ unique up to reordering
such that $$f(\lambda,\mu) = \prod^m_{i=1} (\mu_i \lambda -
\lambda_i \mu)$$ for suitable representatives of homogeneous
coordinates $\lambda_i,\mu_i$. Then by setting $\omega_1 =
\sqrt[m]{\frac{1}{m!}}(\mu_1 e^1 \wedge e^2 + \cdots + \mu_m
e^{2m-1}\wedge e^{2m})$ and $\omega_2 = -
\sqrt[m]{\frac{1}{m!}}(\lambda_1 e^1 \wedge e^2 + \cdots +
\lambda_m e^{2m-1}\wedge e^{2m})$, we have
\begin{align*}
(\lambda \omega_1 + \mu \omega_2)^m &= \frac{1}{m!} [(\mu_1
\lambda - \lambda_1 \mu) e^1 \wedge e^2 + \cdots + (\mu_m \lambda
- \lambda_m \mu) e^{2m-1} \wedge e^{2m}]^m\\
&= \prod^m_{i=1} (\mu_i \lambda - \lambda_i \mu) \, e^1 \wedge
\cdots \wedge e^{2m}\\
&= f(\lambda,\mu) \, e^1 \wedge \cdots \wedge e^{2m}.
\end{align*}
Thus $f_{(\omega_1,\omega_2)} = f$ and $\Phi$ is surjective.
\end{proof}

\begin{corollary} \label{cor:12.2}
$U_1$ consists of infinitely many orbits for $m \geq 4$.
\end{corollary}
\begin{proof}
Note that $$S^m\mathbb{C}^2 /\!/ GL_2 = \left(S^m\mathbb{C}^2-\{0\}\right)/ GL_2 = \overline{\mathcal{M}_{0,m}}$$ where $\mathcal{M}_{0,m}$ is the moduli space of $m$-points in $\mathbb{P}^1$, which is infinite iff $m\geq 4$ according to the fact that $PGL_1$ acts $3$-transitively on $\mathbb{P}^1$. By Proposition \ref{prop:12.1}, $\mathbb{C}^2 \otimes \Lambda^2\mathbb{C}^{2m}/ GL_2 \times SL_{2m}$ is infinite as $\Phi$ is surjective. The result then follows from the one-to-one correspondence between $U_1/ GL_2 \times SL_{2m+1}$ and $\mathbb{C}^2 \otimes \Lambda^2\mathbb{C}^{2m}/ GL_2 \times SL_{2m}$ established above.
\end{proof}

In particular, it forces that the open orbit of $\mathbb{C}^2 \otimes \Lambda^2\mathbb{C}^{2m+1}$ lies in $U_2$, and with a little bit more effort, we see that $U_2$ actually forms a single orbit. The reason is that under the action of $SL_2$ every element $(1,0)\otimes \omega_1 + (0,1)\otimes \omega_2 \in U_2$ can be conjugated so that $\omega_1$, $\omega_2$ are of rank $m$, and that all those full rank elements are inside the same $GL_2 \times SL_{2m+1}$-orbit.

In fact, outside of $\Phi^{-1}(0)$, $$\Phi:\left(\mathbb{C}^2
\otimes \Lambda^2\mathbb{C}^{2m} - \Phi^{-1}(0)\right)/ GL_2
\times SL_{2m} \rightarrow \left(S^m\mathbb{C}^2 -
\{0\}\right)/GL_2$$ is a $m:1$ branched cover of projective
varieties. In particular, when $m \leq 3$, there are only finite
number of orbits upstairs outside the central fibre
$\Phi^{-1}(0)$, while $\Phi^{-1}(0)$ can be identified with
$\left(\mathbb{C}^2 \otimes \Lambda^2\mathbb{C}^{2m-1}\right)/GL_2
\times SL_{2m-1}$. So by backward induction, we see that
$\mathbb{C}^2 \otimes \Lambda^2\mathbb{C}^{2m-1}$ has a finite
number of $GL_2 \times SL_{2m-1}$-orbits for $m \leq 3$. The
result is summerized in the following theorem.

\begin{theorem} \label{thm:12.3}
$(\mathbb{C}^2 \otimes
\Lambda^2\mathbb{C}^{2m+1}, GL_2 \times SL_{2m+1})$ has an open orbit for all $m \geq 1$, and it consists of finite number of orbits if and only if $m \geq 4$.
\end{theorem}

\newpage
\appendix

\section*{\huge Appendix}
\section{$\mathbb{Z}$-Gradations of Semisimple Lie Algebras} \label{sect:3}
This section is devoted to the generalities of $\mathbb{Z}_m$-gradations of semisimple Lie algebras which were encountered when we discussed augmentation of Dynkin diagrams. For the sake of completeness, we have included the proofs of some standard results which can be found in \cite{ref15}.

\begin{definition}
Let $\mathfrak{g}$ be a Lie algebra and $m\in \mathbb{Z}_{\geq 0}$. A $\mathbb{Z}_m$-gradation of $\mathfrak{g}$ is a direct sum decomposition $$\mathfrak{g}= \bigoplus_{i \in \mathbb{Z}_m} \mathfrak{g}_i$$ of $\mathfrak{g}$ into vector subspaces $\mathfrak{g}_i\, (i \in \mathbb{Z}_m)$ satisfying $[\mathfrak{g}_i,\mathfrak{g}_j] \subset \mathfrak{g}_{i+j}$ for all $i,j \in \mathbb{Z}_m$.
\end{definition}

Given a $\mathbb{Z}_m$-gradation $\displaystyle{\mathfrak{g}= \bigoplus_{i \in \mathbb{Z}_m} \mathfrak{g}_i}$ of $\mathfrak{g}$, we see that $\mathfrak{g}_0$ is a subalgebra of $\mathfrak{g}$ and that every $\mathfrak{g}_k(k \neq 0)$ is a $\mathfrak{g}_0$ representation through the adjoint action. Especially, when $\mathfrak{g}$ is complex semisimple and $m=0$\,(i.e.\,$\mathbb{Z}_m=\mathbb{Z}$), there is a more detailed description about the $\mathbb{Z}$-gradation.

\begin{proposition} \label{prop:3.1}
Let $\mathfrak{g}$ be a complex semisimple Lie algebra with a $\mathbb{Z}$-gradation $\displaystyle{\bigoplus_{i\in \mathbb{Z}}\mathfrak{g}_i}$, and let $\kappa: \mathfrak{g} \times \mathfrak{g} \rightarrow \mathbb{C}$ denote the Killing form of $\mathfrak{g}$. Then
\begin{enumerate}
    \item $\kappa(\mathfrak{g}_i,\mathfrak{g}_j)=0$ whenever $i+j \neq 0$.
    \item $\kappa|_{\mathfrak{g}_i \times \mathfrak{g}_{-i}}$ is nondegenerate for all $i \in \mathbb{Z}$; in particular, it implies that $\mathfrak{g}_0$ is a reductive subalgebra of $\mathfrak{g}$.
\end{enumerate}
\end{proposition}
\begin{proof}
\begin{enumerate}
    \item Pick any $X \in \mathfrak{g}_i, Y \in \mathfrak{g}_j$, then for every $k \in \mathbb{Z}$ $$(\text{ad}_X \circ \, \text{ad}_Y)^m(\mathfrak{g}_k) \subset \mathfrak{g}_{k+ m(i+j)}.$$ Since $i+j\neq 0$, for sufficiently large $m$, $\mathfrak{g}_{k+ m(i+j)}=0$. Hence $\text{ad}_X \circ \,\text{ad}_Y$ is nilpotent and $\kappa(X,Y)= \text{ Tr}(\text{ad}_X \circ \,\text{ad}_Y)=0$.
    \item For every nonzero $X \in \mathfrak{g}_i$, there exists an $Y\in \mathfrak{g}$ such that $\kappa(X,Y)\neq 0$ as $\kappa$ is nondegenrate on $\mathfrak{g}$. Now let $Y_{-j}$ be the component of $Y$ in $\mathfrak{g}_j$ for $j \in \mathbb{Z}$. Then in view of (i), $$\kappa(X,Y)= \sum_{j\in \mathbb{Z}}\kappa(X,Y_j) = \kappa(X,Y_{-i}) \neq 0.$$ Thus $\kappa|_{\mathfrak{g}_i \times \mathfrak{g}_{-i}}$ is nondegenerate.
\end{enumerate}
\end{proof}
\section{Basic Facts about Algebraic Groups} \label{sect:4}

Let $\mathcal{V}$ be a complex $G$-variety, i.e. a complex algebraic variety with a continuous group homomorphism $\pi$ from $G$ to the group of biregular morphisms $\mathbb{C}[\mathcal{V}]^*$ on $\mathcal{V}$. For each $x\in \mathcal{V}$, we can form the orbit $G \cdot x$ and consider the orbit map $$\pi_x : G \rightarrow G \cdot x.$$

\begin{proposition} \label{prop:4.1}
Let $G, \mathcal{V}, \pi, \pi_x$ be defined as above.
\begin{enumerate}
    \item For each $x \in \mathcal{V}$, the orbit closure $\overline{G \cdot x}$ is a subvariety of $\mathcal{V}$. Moreover, if $\mathcal{V}$ is affine, then so is $\overline{G \cdot x}$.
    \item $G \cdot x$ is open in $\overline{G \cdot x}$. In particular, there is a natural structure of smooth algebraic variety on $G \cdot x$.
    \item The orbit map $\pi_x$ is a surjective morphism of varieties.
\end{enumerate}
\end{proposition}

Note that the identity component $G^0$ of $G$ is connected, which is equivalent to $G^0$ being
irreducible. We have, for every $x\in \mathcal{V}$, $G^0 \cdot x=
\pi_x(G^0)$ is irreducible, as $\pi_x$ is a surjective morphism. In
general, we can write $G \cdot x$ as a finite union of
$G^0$-orbits in $G$; these $G^0$-orbits are both connected and
irreducible components of $G\cdot x$.

\begin{proposition} \label{prop:4.2}
Let $G_x:=\{g\in G| g \cdot x=x\}$ denote the stabilizer (also
called isotropy subgroup) of $x \in \mathcal{V}$. Then $G_x$ is a
closed subgroup of $G$ and $\pi_x$ induces an isomorphism
$$\overline{\pi_x}: G/G_x \longrightarrow G \cdot x.$$ In
particular, we have $$ \dim G\cdot x = \dim G - \dim G_x.$$
\end{proposition}

\begin{corollary}
For every $x \in \mathcal{V}$, $G\cdot x$ is a smooth
equidimensional algebraic variety. More precisely, all irreducible
components of $G \cdot x$ are smooth subvarieties having the same
dimension $\dim G - \dim G\cdot x$.
\end{corollary}

Note that $\overline{G \cdot x}$ is clearly stable under the
action of $G$, hence it is a union of $G$-orbits. In
particular, this enables us to define a partial ordering on the
set of $G$-orbits.

\begin{definition}
For any pair of elements $x,y \in \mathcal{V}$, we say that $G
\cdot y$ is less than $G \cdot x$, denoted by $G \cdot y \prec G
\cdot x$, if $G \cdot y \subseteq \overline{G \cdot x}$. This
yields are partial ordering, called the closure ordering, on
the set $G \backslash \mathcal{V}$ of $G$-orbits in $\mathcal{V}$.
\end{definition}


Now let $\mathcal{V}$ be defined over $\mathbb{R}$, so that the set of $\mathbb{R}$-rational points $\mathcal{V}_\mathbb{R}$ is a variety over $\mathbb{R}$. We have the following fundamental result of Whitney \cite{ref16}:

\begin{theorem} \label{thm:4.4}
Let $\mathcal{V}$ be a complex algebraic variety defined over $\mathbb{R}$. Then the set of $\mathbb{R}$-rational points $\mathcal{V}_\mathbb{R}$ of $\mathcal{V}$ decomposes into a finite number of connected components.
\end{theorem}

\begin{corollary}
Let $G$ be a connected complex algebraic group defined over $\mathbb{R}$. Then $G_\mathbb{R}$ has a finite number of connected components.
\end{corollary}

\begin{theorem} \label{thm:4.6}
Let $G$ be a complex reductive algebraic group defined over $\mathbb{R}$ and $V$ be a representation of $G$ with finite number of $G$-orbits whose restriction to $G_\mathbb{R}$ is of real type. Then the real representation $V_\mathbb{R}$ of $G_\mathbb{R}$ has a finite number of $G_\mathbb{R}$-orbits. In particular, $V_\mathbb{R}$ has an open $G_\mathbb{R}$-orbit.
\end{theorem}
\begin{proof}
Note that every $G$-orbit in $V$ is also stable under $G_\mathbb{R}$. So it suffices to show that every $G$-orbit intersects $V_\mathbb{R}$ with a finite number of $G_\mathbb{R}$-orbits. Now fix any $G$-orbit $\mathcal{O}$ in $V$. Then for every $v \in \mathcal{O} \cap V_\mathbb{R}$, $$T_v(\mathcal{O} \cap V_\mathbb{R}) \subset \mathfrak{g}\cdot v \cap V_\mathbb{R} = \mathfrak{g}_{\mathbb{R}} = T_v(G_\mathbb{R} \cdot w).$$ On the other hand, we have $\mathcal{O} \cap V_\mathbb{R} \supset G_\mathbb{R} \cdot v$. Thus $T_v(\mathcal{O} \cap V_\mathbb{R})= T_v(G_\mathbb{R} \cdot v)$. It follows that $G_\mathbb{R}\cdot v$ is open in $\mathcal{O} \cap V_\mathbb{R}$ and is also closed in $\mathcal{O} \cap V_\mathbb{R}$ since its complement is union of such $G_\mathbb{R}$-orbits. Hence we conclude that $G_\mathbb{R} \cdot v$ is a union of connected components of $\mathcal{O} \cap V_\mathbb{R}$. Now by Theorem \ref{thm:4.4}, number of connected components of $\mathcal{O} \cap V_\mathbb{R}$ must be finite. As the number of connected components of $\mathcal{O} \cap V_\mathbb{R}$ must exceed the number of $G_\mathbb{R}$-orbits in $\mathcal{O} \cap V_\mathbb{R}$, there can only have finite number of $G_\mathbb{R}$-orbits in $\mathcal{O} \cap V_\mathbb{R}$.
\end{proof}

\newpage
\section{Tables} \label{sect:13}
\begin{table}[h]
    \begin{center}
  \caption{Table for irreducible PVS of parabolic type.}\label{table:2}
  \begin{tabular}{|c|c|c|c|}
    \hline
    $G$ & $\alpha_0$ & \phantom{$\displaystyle{\bigoplus}$}$(G_0,V)$\phantom{$\displaystyle{\bigoplus}$} & $\nu(\mathfrak{g},V)$\\
    \hline
    $A_n$ & & $\xymatrix{\circ \ar@{-}[r]_{\!\!\!\!\!\!\!\!\!\!\!\!\!\!\!\!\!\!\!1} &\circ \ar@{-}[r]_{\!\!\!\!\!\!\!\!\!\!\!\!\!\!\!\!\!\!\!2} &\cdots \ar@{-}[r] &\circ \ar@{-}[r]_{\!\!\!\!\!\!\!\!\!\!\!\!\!\!\!\!\!\!\!n-1} &\circ \ar@{-}[l]^{\,\,\,\,\qquad \qquad n\phantom{-1}}}$ & \phantom{$\displaystyle{\bigoplus}_{\mathbb{Q}}^\mathbb{Q}$}\\
    \phantom{$\displaystyle{\bigoplus}$} & $1$ & $(GL_n, \mathbb{C}^n)$ & $(1)$\\
    \phantom{$\displaystyle{\bigoplus}$} & $k$ & $(GL_k \times SL_{n-k}, \mathbb{C}^k \otimes (\mathbb{C}^{n-k})^*)$ \begin{tiny}$2 \leq k \leq \frac{n+1}{2}$\end{tiny} & $(1,1)$\\
    \hline
    \small{$B_n$} & & $\xymatrix{\circ \ar@{-}[r]_{\!\!\!\!\!\!\!\!\!\!\!\!\!\!\!\!\!\!\!1} & \circ \ar@{-}[r]_{\!\!\!\!\!\!\!\!\!\!\!\!\!\!\!\!\!\!\!\!2} & \cdots \ar@{-}[r] &\circ \ar@2{->}[r]_{\,\,\,\,\,\,\,n-1\qquad n} &\circ\\}$ & \phantom{$\displaystyle{\bigoplus}_{\mathbb{Q}}^\mathbb{Q}$}\\
    \phantom{$\displaystyle{\bigoplus}$} & $1$ & $(GL_1 \times SO_{2n-1}, \mathbb{C} \otimes \mathbb{C}^{2n-1})$ & (1)\\
    \phantom{$\displaystyle{\bigoplus}$} & $k$ & $(GL_k \times SO_{2n-2k+1}, \mathbb{C}^k \otimes \mathbb{C}^{2n-2k+1})$ \begin{tiny}$2 \!\leq\! k \!\leq \!n\!-\!\!1$\end{tiny} & $(1,1)$\\
    \phantom{$\displaystyle{\bigoplus}$} & $n$ & $(GL_n, \mathbb{C}^n)$ & $(2)$\\
    \hline
    $C_n$ & & $\xymatrix{\circ \ar@{-}[r]_{\!\!\!\!\!\!\!\!\!\!\!\!\!\!\!\!\!\!\!\!1} & \circ \ar@{-}[r]_{\!\!\!\!\!\!\!\!\!\!\!\!\!\!\!\!\!\!\!\!2} & \cdots \ar@{-}[r]_{\,\,\,\,\,\qquad \qquad n-1} &\circ &\circ \ar@2{->}[l]^{\,\,\,\,\,\,\,\qquad \qquad n\phantom{-1}}\\}$ & \phantom{$\displaystyle{\bigoplus}_{\mathbb{Q}}^\mathbb{Q}$}\\
    \phantom{$\displaystyle{\bigoplus}$} & $1$ & $(GL_1 \times Sp_{n-1}, \mathbb{C} \otimes \mathbb{C}^{2n-2})$ & $(1)$\\
    \phantom{$\displaystyle{\bigoplus}$} & $k$ & $(GL_k \times Sp_{n-k}, \mathbb{C}^k \otimes \mathbb{C}^{2n-2k})$ \begin{tiny}$2 \leq k \leq n-2$\end{tiny} & $(1,1)$\\
    \phantom{$\displaystyle{\bigoplus}$} & $n-1$ & $(GL_{n-1} \times SL_2, \mathbb{C}^n-1 \otimes \mathbb{C}^2)$ & $(1,2)$\\
    \phantom{$\displaystyle{\bigoplus}$} & $n$ & $(GL_n, S^2\mathbb{C}^n)$ & $(1)$\\
    \hline
    $D_n$ & & $\xymatrix{\circ \ar@{-}[r]_{\!\!\!\!\!\!\!\!\!\!\!\!\!\!\!\!\!\!\!1} &\circ \ar@{-}[r]_{\!\!\!\!\!\!\!\!\!\!\!\!\!\!\!\!\!\!\!\!2} &\cdots \ar@{-}[r] &\circ \ar@{-}[r]_{\!\!\!\!\!\!\!\!\!\!\!n-2} &\circ \ar@{-}[l]^{\!\!\qquad\qquad \qquad n-1}\\
&&&\circ \ar@{-}[u] \ar@{-}_{n\phantom{-1}}\\}$ & \phantom{$\displaystyle{\bigoplus}$}\\
    \phantom{$\displaystyle{\bigoplus}$} & $1$ & $(GL_1 \times SO_{2n-2}, \mathbb{C} \otimes \mathbb{C}^{2n-2})$ & $(1)$\\
    \phantom{$\displaystyle{\bigoplus}$} & $k$ & $(GL_k \times SO_{2n-2k}, \mathbb{C}^k \otimes \mathbb{C}^{2n-2k})$ \begin{tiny}$2 \leq k \leq n-2$\end{tiny} & $(1,1)$\\
    \phantom{$\displaystyle{\bigoplus}$} & $n$ & $(GL_n, \Lambda^2\mathbb{C}^n)$ & $(1)$\\
    \hline
    \end{tabular}
    \end{center}
\end{table}
\begin{table*}
    \begin{center}
    \begin{tabular}{|c|c|c|c|}
    \hline
     $E_6$ & & $\xymatrix{\circ \ar@{-}[r]_{\!\!\!\!\!\!\!\!\!\!\!\!\!\!\!\!\!\!\!1} &\circ \ar@{-}[r]_{\!\!\!\!\!\!\!\!\!\!\!\!\!\!\!\!\!\!\!\!3} &\circ \ar@{-}[r]_{\!\!\!\!\!\!\!\!\!\!\!\!\!\!\!\!\!\!\!\!4} &\circ \ar@{-}[r]_{\!\!\!\!\!\!\!\!\!\!\!\!\!\!\!\!\!\!\!\!5} &\circ \ar@{-}[l]^{\,\,\,\qquad \qquad 6}\\
&&\circ \ar@{-}[u] \ar@{-}_{2}\\}$ & \phantom{$\displaystyle{\bigoplus}$}\\
    \phantom{$\displaystyle{\bigoplus}$} & $1$ & $(GL_1 \times Spin_{10},\mathbb{C} \otimes S^+)$ & $(1)$\\
    \phantom{$\displaystyle{\bigoplus}$} & $3$ & $(GL_2 \times SL_5, \mathbb{C}^2 \otimes \Lambda^2\mathbb{C}^5)$ & $(1,1)$\\
    \phantom{$\displaystyle{\bigoplus}$} & $4$ & $(GL_2 \times SL_3^2, \mathbb{C}^2 \otimes \mathbb{C}^3 \otimes \mathbb{C}^3)$ & $(1,1,1)$\\
    \phantom{$\displaystyle{\bigoplus}$} & $2$ & $(GL_6, \Lambda^3 \mathbb{C}^6)$ & $(1)$\\
    \hline
    $E_7$ & & $\xymatrix{\circ \ar@{-}[r]_{\!\!\!\!\!\!\!\!\!\!\!\!\!\!\!\!\!\!\!\!1} &\circ \ar@{-}[r]_{\!\!\!\!\!\!\!\!\!\!\!\!\!\!\!\!\!\!\!\!3} &\circ \ar@{-}[r]_{\!\!\!\!\!\!\!\!\!\!\!\!\!\!\!\!\!\!\!\!4} &\circ \ar@{-}[r]_{\!\!\!\!\!\!\!\!\!\!\!\!\!\!\!\!\!\!\!\!5} &\circ \ar@{-}[r]_{\!\!\!\!\!\!\!\!\!\!\!\!\!\!\!\!\!\!\!\!6} &\circ \ar@{-}[l]^{\,\,\,\qquad \qquad 7}\\
&&\circ \ar@{-}[u] \ar@{-}_{2}\\}$ & \phantom{$\displaystyle{\bigoplus}$}\\
     & $1$ & $(GL_1 \times Spin_{12}, \mathbb{C} \otimes S^+)$ & $(1)$\\
     & $3$ & $(GL_2 \times SL_6, \mathbb{C}^2 \otimes \Lambda^2 \mathbb{C}^6)$ & $(1,1)$\\
     & $4$ & $(GL_2 \times SL_3 \times SL_4, \mathbb{C}^2 \otimes \mathbb{C}^3 \otimes \mathbb{C}^4)$ & $(1,1,1)$\\
     & $5$ & $(GL_3 \times SL_5, \mathbb{C}^3 \otimes \Lambda^2 \mathbb{C}^5)$ & $(1,1)$\\
     & $6$ & $(GL_2 \times Spin_{10}, \mathbb{C}^2 \otimes S^+)$ & $(1,1)$\\
     & $7$ & $(GL_1 \times E_6, \mathbb{C} \otimes \mathbb{C}^{27})$ & $(1)$\\
     & $2$ & $(GL_7, \Lambda^3\mathbb{C}^7)$ & $(1)$\\
    \hline
    $E_8$ & & $\xymatrix{\circ \ar@{-}[r]_{\!\!\!\!\!\!\!\!\!\!\!\!\!\!\!\!\!\!\!\!1} &\circ \ar@{-}[r]_{\!\!\!\!\!\!\!\!\!\!\!\!\!\!\!\!\!\!\!\!3} &\circ \ar@{-}[r]_{\!\!\!\!\!\!\!\!\!\!\!\!\!\!\!\!\!\!\!\!4} &\circ \ar@{-}[r]_{\!\!\!\!\!\!\!\!\!\!\!\!\!\!\!\!\!\!\!\!5} &\circ \ar@{-}[r]_{\!\!\!\!\!\!\!\!\!\!\!\!\!\!\!\!\!\!\!\!6} &\circ \ar@{-}[r]_{\!\!\!\!\!\!\!\!\!\!\!\!\!\!\!\!\!\!\!\!7} &\circ \ar@{-}[l]^{\,\,\,\qquad \qquad 8}\\
&&\circ \ar@{-}[u] \ar@{-}_{2}\\}$ & \phantom{$\displaystyle{\bigoplus}$}\\
     & $1$ & $(GL_1 \times Spin_{14}, \mathbb{C} \otimes S^+)$ & $(1)$\\
     & $3$ & $(GL_2 \times SL_7, \mathbb{C}^2 \otimes \Lambda^2\mathbb{C}^7)$ & $(1,1)$\\
     & $4$ & $(GL_2 \times SL_3 \times SL_5, \mathbb{C}^2 \otimes \mathbb{C}^3 \otimes \mathbb{C}^5)$ & $(1,1,1)$\\
     & $5$ & $(GL_4 \times SL_5, \mathbb{C}^4 \otimes \Lambda^2\mathbb{C}^5)$ & $(1,1)$\\
     & $6$ & $(GL_3 \times Spin_{10}, \mathbb{C}^3 \otimes S^+)$ & $(1,1)$\\
     & $7$ & $(GL_2 \times E_6, \mathbb{C}^2 \otimes \mathbb{C}^{27})$ & $(1,1)$\\
     & $8$ & $(GL_1 \times E_7, \mathbb{C} \otimes \mathbb{C}^{56})$ & $(1)$\\
     & $2$ & $(GL_8, \Lambda^3 \mathbb{C}^8)$ & $(1)$\\
    \hline
    $F_4$ & & $\xymatrix{\circ \ar@{-}[r]_{\!\!\!\!\!\!\!\!\!\!\!\!\!\!\!\!\!\!\!\!1} &\circ \ar@2{->}[r]_{\!\!\!\!\!\!\!\!\!\!\!\!\!\!\!\!\!\!\!\!\!\!2} &\circ \ar@{-}[r]_{\!\!\!\!\!\!\!\!\!\!\!\!\!\!\!\!\!\!\!\!3} &\circ \ar@{-}[l]^{\,\,\,\,\qquad \qquad 4}\\}$ & \phantom{$\displaystyle{\bigoplus}_{\mathbb{Q}}^\mathbb{Q}$}\\
     & $4$ & $(GL_1 \times Spin_7, \mathbb{C} \otimes S)$ & $(1)$\\
     & $3$ & $(GL_2 \times SL_3, \mathbb{C}^2 \otimes \mathbb{C}^3)$ & $(1,2)$\\
     & $2$ & $(GL_2 \times SL_3, \mathbb{C}^2 \otimes S^2 \mathbb{C}^3)$ & $(1,1)$\\
     & $1$ & $(GL_1 \times Sp_3, \mathbb{C} \otimes \Lambda^3_0 \mathbb{C}^6)$ & $(1)$\\
    \hline
    $G_2$ & & $\xymatrix{\circ \ar@3{<-}[r]_{\!\!\!\!\!\!\!\!\!\!\!\!\!\!\!\!\!\!\!\!1} &\circ \ar@{-}_{\,\,\,\,2}}$ & \phantom{$\displaystyle{\bigoplus}_{\mathbb{Q}}^\mathbb{Q}$}\\
     & $1$ & $(GL_2, \mathbb{C}^2)$ & $(3)$\\
     & $2$ & $(GL_2, S^3\mathbb{C}^2)$ & $(1)$\\
    \hline
  \end{tabular}
  \end{center}
\end{table*}

\begin{table}
  \centering
  \caption{Table for twisted affine type.}\label{table:3}
  \begin{tabular}{|c|c|c|c|}
    \hline
    \phantom{$\displaystyle{\bigoplus}$}$\mathfrak{g}$\phantom{$\displaystyle{\bigoplus}$} & \phantom{$\displaystyle{\bigoplus}$}$(G_0,V)$\phantom{$\displaystyle{\bigoplus}$}\\
    \hline
    \phantom{$\displaystyle{\bigoplus}$}$\mathfrak{so}_{10}$\phantom{$\displaystyle{\bigoplus}$} & $(GL_1 \times G_2, \mathbb{C} \otimes \mathbb{C}^7)$ \\
    \phantom{$\displaystyle{\bigoplus}$}$E_6$\phantom{$\displaystyle{\bigoplus}$} & $(GL_2 \times G_2, \mathbb{C}^2 \otimes \mathbb{C}^7)$ \\
    \phantom{$\displaystyle{\bigoplus}$}$\mathfrak{so}_{12}$\phantom{$\displaystyle{\bigoplus}$} & \phantom{aa}$(GL_2 \times Spin_7, \mathbb{C}^2 \otimes S)$\phantom{aa}\\
    \phantom{$\displaystyle{\bigoplus}$}$E_7$\phantom{$\displaystyle{\bigoplus}$} & $(GL_3 \times Spin_7, \mathbb{C}^3 \otimes S)$ \\
    \phantom{$\displaystyle{\bigoplus}$}$E_6$\phantom{$\displaystyle{\bigoplus}$} & $(GL_1 \times Spin_9, \mathbb{C} \otimes S)$ \\
    \phantom{$\displaystyle{\bigoplus}$}$E_7$\phantom{$\displaystyle{\bigoplus}$} & $(GL_1 \times Spin_{11}, \mathbb{C} \otimes S)$ \\
    \hline
  \end{tabular}
\end{table}

\newpage
\noindent \textbf{Acknowledgements.} Both authors thank H.
Rubenthaler and J.M.Landsberg for informing us that many of our
results were previously known, and pointing out some minor
mistakes in the preliminary version of our paper. The second
author is partially supported by a RGC research grant from the
Hong Kong Government.

\end{document}